\documentclass[a4paper,12pt]{amsart}

\usepackage[
  margin=30mm,
  marginparwidth=25mm,     
  marginparsep=2mm,       
  bottom=25mm,
  ]{geometry}

\usepackage[T1]{fontenc}
\usepackage[utf8]{inputenc}
\usepackage[english]{babel}
\usepackage{amsmath,amssymb,amsthm,mathtools}

\expandafter\let\expandafter\xbf\csname bfseries \endcsname
\expandafter\let\expandafter\xmd\csname mdseries \endcsname
\let\xbar\bar
\usepackage{allrunes}
\let\bar\xbar
\expandafter\let\csname bfseries \endcsname\xbf
\expandafter\let\csname mdseries \endcsname\xmd

\usepackage{latexsym}
\usepackage{delarray}
\usepackage{bbm}
\usepackage{scrextend}
\usepackage{hyperref}
\usepackage[pdftex,usenames,dvipsnames]{xcolor}
\usepackage{datetime}
\usepackage{mathrsfs}
\usepackage{enumitem}  
\usepackage{tikz}
\usetikzlibrary{fadings}

\newtheorem{Th}{Theorem}[section]
\newtheorem{Prop}[Th]{Proposition}
\newtheorem{Lem}[Th]{Lemma}
\newtheorem{Cor}[Th]{Corollary}
\newtheorem{Rem}[Th]{Remark}
\newtheorem{Def}[Th]{Definition}

\newcommand{\eps}{\varepsilon}

\newcommand{\Rl}{\mathbb{R}}

\newcommand{\Rn}{\mathbb{R}^{n}}

\renewcommand{\d}{\partial}

\def\bra#1{{\langle{#1}\rangle}}

\newcommand{\jap}[1]{\left\langle #1\right\rangle}
\newcommand{\bmo}{\mathrm{bmo}}
\newcommand{\phase}{\varphi}

\newcommand{\dd}{\, \mathrm{d}}
\newcommand{\ddd}{\,\text{\rm{\mbox{\dj}}}}

\renewcommand{\SS}{{\mathscr{S}}}

\newcommand{\at}{\mathfrak{a}}

\newcommand{\brkt}[1]{\Big({#1}\Big)}
\newcommand{\set}[1]{\left\{{#1}\right\}}
\newcommand{\norm}[1]{\Big\|#1\Big\|}
\newcommand{\abs}[1]{\Big |#1\Big |}
\newcommand{\diffcases}[1]{\begin{cases}#1\end{cases}}

\newcommand{\eq}[1]{
    \begin{align*}
        #1
    \end{align*}
}

\newcommand{\nm}[2]{\begin{equation}\label{#1}
\begin{split}
#2
\end{split}
\end{equation}}

\DeclareMathOperator{\supp}{supp}

\allowdisplaybreaks 

\title[Oscillatory integral operators]
{Regularity of oscillatory integral operators}

\author[A. Israelsson]{Anders Israelsson}
\author[T. Mattsson]{Tobias Mattsson}
\author[W. Staubach]{Wolfgang Staubach}

\address{\newline
      Anders Israelsson, Tobias Mattsson,
       Wolfgang Staubach \newline
       Department of  Mathematics, Uppsala University, \newline
       S-751 06 Uppsala, Sweden}
       \email{anders.israelsson@math.uu.se, tobias.mattsson@math.uu.se, wulf@math.uu.se}

 \thanks{
 The second author is supported by the Knut and Alice Wallenberg Foundation.}



\begin{document}


\maketitle

\begin{abstract}
In this paper, we establish the global boundedness of oscillatory integral operators on Besov-Lipschitz and Triebel-Lizorkin spaces, with amplitudes in general $S^m_{\rho,\delta}(\Rn)$-classes and non-degenerate phase functions in the class $\textart F^k$. Our results hold for a wide range of parameters $0\leq\rho\leq1$, $0\leq\delta<1$, $0<p\leq\infty$, $0<q\leq\infty$ and $k>0$. We also provide a sufficient condition for the boundedness of operators with amplitudes in the forbidden class $S^m_{1,1}(\Rn)$ in Triebel-Lizorkin spaces.
\end{abstract}

\hspace*{0.75cm}\\
\tableofcontents

\tableofcontents

\section{Introduction}

This paper is devoted to the investigation of the global regularity of oscillatory integral operators (here referred to as OIOs) of the form
\begin{equation*}
	T_a^\varphi f(x) = \frac{1}{(2\pi)^n} \int_{\Rl^n} e^{i\varphi(x,\xi)}\,a(x,\xi)\,\widehat f (\xi) \dd\xi,
\end{equation*}
with amplitudes in the general H\"ormander class $S^m_{\rho, \delta}(\Rl^n)$ (Definition \ref{symbol class Sm}), 
on Besov-Lipschitz $B^{s}_{p,q}(\Rl^n)$ and Triebel-Lizorkin $F^{s}_{p,q}(\Rl^n)$ spaces of order $s\in \Rl$ with $0<p\leq\infty$ and $0< q\leq \infty.$ Throughout the paper, we are assuming that the phase function $\varphi(x, \xi)$ is in the class $\textart F^k$ for some $k>0$ (Definition \ref{def:fk}) and is strongly non-degenerate (Definition \ref{nondeg phase}). Here we note that we are not assuming any homogeneity in the $\xi$ variable of the phase function $\varphi$, as is the case for Fourier integral operators where $\varphi(x,\xi)$ is assumed to be positively homogeneous of degree one in $\xi$.\\

The motivation for the study of the OIOs that are investigated in this paper comes from the theory of partial differential equations where
phase functions $\varphi(x,\xi)= x\cdot \xi+ \phi(\xi)$ frequently appear in the study of dispersive equations. Indeed $\phi(\xi)= |\xi|^{1/2}$ corresponds to the water-wave equation, $\phi(\xi)= |\xi|^{2}$ corresponds to the Schr\"odinger equation,
while $\phi(\xi)= |\xi|^{3}$ and $\phi(\xi)= \xi|\xi|$ (both in dimension one) corresponds to Airy and Benjamin-Ono equations respectively.

The results that are obtained in this paper also accommodate for instance the case of variable coefficient Schr\"odinger equations through the earlier investigations of B. Helffer and  D. Robert \cite{Helf-Rob}, Helffer \cite{Helffer} and the work of E. Cordero F. Nicola and L. Rodino   \cite{CNR}, \cite{CoNiRo}.\\

In \cite{CISY}, the second and the third author in collaboration with A. Castro and M. Yerlanov established an $L^p-L^q$, ($1< p\leq q<\infty$) regularity theory for OIOs in the following two cases: (1) when the amplitude is in $S^m_{1,0}(\Rn)$ and the phase function is in $\textart F^k$ with $0<k<1$. (2) when the amplitude is in $S^m_{0,0}(\Rn)$ and the phase function is in $\textart F^k$ with $k\geq 1$. The authors of \cite{CISY} also went beyond the scope of $L^p$-spaces and investigated the regularity of OIOs in classical function spaces such as Besov-Lipschitz and Triebel-Lizorkin spaces.\\

\noindent In \cite{PRS} M. Pramanik, K. Rogers and A. Seeger  (see Theorem \ref{thm:PRS_proto} below), proved a Calder\'on-Zygmund-type estimate with far-reaching applications, including the regularity of Radon transforms and Fourier integral operators. In that paper, the authors also considered local Fourier integral operators $T_a^{\varphi}$ where $a(x,\xi)\in S^{m}_{1,0}(\Rn)$ is compactly supported in $x$ and $\varphi(x, \xi)$ is non-degenerate on the support of $a(x,\xi).$ Using their Calder\'on-Zygmund estimate in \cite{PRS}, they showed that if \(n \geq 2,\,2<p<\infty,\, q>0\), then \(T_a^{\varphi}: F^{0}_{p,p}(\mathbb{R}^{n}) \rightarrow F^{0}_{p,q}(\mathbb{R}^{n})\), provided that $m=-(n-1)\big|\frac{1}{p}-\frac{1}{2}\big|.$\\

In this paper we prove a variation of Theorem \ref{thm:PRS_proto} (see Theorem \ref{thm:PRS}), which is on the one hand more suitable for extensions to the case of quasi-Banach scales of Besov-Lipschitz and Triebel-Lizorkin spaces, and on the other also fits well for the applications to the regularity theory of OIOs. More specifically using certain decomposition in the frequency space and rather intricate kernel estimates for oscillatory integral operators, a composition theorem for the action of parameter-dependent pseudodifferential operators on OIOs, atomic and molecular decompositions of Triebel-Lizorkin spaces in the spirit of Frazier-Jawerth \cite{FJ:phi-transform, FJ:discrete-transform, FJ-phi-and-wavelet-transform}, and vector-valued inequalities of the type provided in Theorem \ref{thm:PRS}, we manage to get a significant extension of the results  in \cite{CISY}. These extensions are both in terms of the types of oscillatory integral operators, and also the scales of the function spaces on which the operators act. Figure \ref{pic:TLendpointresults} illustrates the extensions that are obtained here.\\ 
More specifically in Theorem \ref{thm:TL_tdphase_k} we obtain the following result.
Let $0\leq \rho\leq 1$ and $0\leq\delta<1$, $0<\mu\leq 1,$ $k>0$, 
    $\varkappa=\min(\rho,1-k)$ and set 
\begin{equation*}
    m(p):= m_\varkappa(p)+ \zeta.
\end{equation*}
where
$$m_\varkappa(p):=-n(1-\varkappa)\Big|\frac{1}{p}-\frac{1}{2}\Big|$$
and
$$\zeta:=\min\Big(0,\frac{n}{2}(\rho-\delta)\Big).$$ Assume that $a\in S^{m(p)}_{\rho, \delta}(\Rl^n)$, and suppose furthermore that $\varphi\in \textart F^k$ {is {SND}} and satisfies the $L^2$-condition of Definition \ref{L2 conditions} and the {LF}$(\mu)$-condition of Definition \eqref{Def:LFmu}. Let $T_a^\varphi$  be the associated {OIO}.

    If $s\in \Rl$, $0<q\leq \infty$, and either one of the following cases holds
    \begin{align*}
    (i)\quad &2<p<\infty \text{  when }0<q\leq p,\\
    (ii) \quad &\frac{n}{n+\mu}<p<2 \text{ when }p\leq q,\\
    (iii)\quad &p=q=2,
    \end{align*}
    then the OIO $T_a^\varphi$ is bounded from $F_{p,q}^{s}(\Rn)$ to $F_{p,q}^s(\Rn).$ Moreover in Theorem \ref{thm:BL_tdphase_k_less_than_one} we also show the boundedness of $T_a^\varphi$ on Besov-Lipschitz spaces $B_{p,q}^s(\Rn),$ for any $\frac{n}{n+\mu}<p<\infty$ and any $0<q\leq\infty.$ These results are currently the most general regularity results for  oscillatory integral operators with amplitudes in a general H\"ormander class, which include the majority of OIOs that appear in the theory of partial differential equations.\\

    The case of operators with amplitudes in the forbidden H\"ormander class \linebreak$S^{-kn|1/p-1/2|}_{1, 1}(\Rl^n)$ is excluded in the results above since these operators do not in general even allow $L^2$-boundedness. However, for $s>n\big(\frac{1}{\min\{1, p, q\}}-1\big)$ we are able to show (Theorem \ref{thm:Sobolev_oio}) the boundedness of OIOs on $F_{p,q}^{s}(\Rn)$ under either of the cases (\emph{i}), (\emph{ii}) and (\emph{iii}) above.\\

The paper is organized as follows; in Section \ref{prelim} we provide the necessary preliminaries from the theory of oscillatory integral operators and the theory of function spaces. Here the reader will also find some of the basic results that are used throughout the paper. In Section \ref{basic kernel estimates}
we prove a basic kernel estimate which lies at the ground for the establishment of Besov-Lipschitz boundedness of OIOs. This is done by utilizing a particular frequency-space decomposition adapted to the OIOs with phase functions in the class $\textart F^k$. In Section \ref{L2 section} we prove the basic global $L^2$-boundedness result for OIOs with strongly non-degenerate phase functions and amplitudes in general H\"ormander classes. This is done by using a continuous version of the almost orthogonality method and Cotlar-Stein's lemma. In Section \ref{sec:hplp} we prove the $h^p\to L^p$ ($p<1)$ and $L^\infty \to \mathrm{bmo}$ boundedness of OIOs. In Sections \ref{subsec:OIO_lowP_transference} and \ref{subsec:OIO_PRS_transference} the Hardy space results are transferred to Triebel-Lizorkin spaces. Section \ref{main TL estim section} contains our main result concerning the boundedness of OIOs on Triebel-Lizorkin spaces. This section includes the results for both classical and forbidden amplitudes. In Section \eqref{main results in BL} we conclude the paper by proving our main results on the regularity of OIOs in Besov-Lipschitz spaces.\\

{{\bf{Acknowldgements.}}
The second author is  supported by the Knut and Alice Wallenberg Foundation.
The authors are also grateful to Andreas Str\"ombergsson for his support and encouragement.}

\section{Preliminaries}\label{prelim}

In this section, we introduce the necessary background and preliminary results that will be required for the development of our main results. These results, while not the focus of our paper, are crucial for the understanding and appreciation of our work.\\

We begin by introducing the key concepts and definitions relating to oscillatory integral operators. Next, we present Theorem \ref{thm:left composition with pseudo} on the composition of pseudodifferential operators and oscillatory integral operators and then basic facts about Triebel-Lizorkin and Besov-Lipschitz spaces which will be used in the proof of our main results.\\

By the end of this section, the reader should have a strong foundation and be well-equipped to move on to the main results of our paper.\\

As is common practice, we will denote positive constants in the inequalities by $C$, which can be determined by known parameters in a given situation but whose
value is not crucial to the problem at hand. Such parameters in this paper would be, for example, $m$, $p$, $s$, $n$,  and the constants connected to the seminorms of various amplitudes or phase functions. The value of $C$ may differ
from line to line, but in each instance could be estimated if necessary. We also write $a\lesssim b$ as shorthand for $a\leq Cb$ and moreover will use the notation $a\sim b$ if $a\lesssim b$ and $b\lesssim a$.\\

\subsection{Basic facts related to oscillatory integral operators}
\quad\\

\noindent Oscillatory integral operators play a central role in a wide range of mathematical fields, including harmonic analysis, partial differential equations, and the study of singular integral operators. In this section, we provide a brief overview of some basic facts related to oscillatory integral operators that will be used throughout the rest of our paper.\\

We begin by introducing the definition of an oscillatory integral operator and discussing some of its basic properties. Next, we present a key lemma and theorem that provide useful information on the compositions of these operators.\\

\begin{Def}\label{symbol class Sm}
Let $m\in \Rl$ and $\rho, \delta \in [0,1]$. An \textit{amplitude} \emph{(}symbol\emph{)} $a(x,\xi)$ in the class $S^m_{\rho,\delta}(\Rl^n)$ is a function $a\in \mathcal{C}^\infty (\Rl^n\times \Rl^n)$ that verifies the estimate
\begin{equation*}
\left|\partial_{\xi}^\alpha \partial_{x}^\beta a(x,\xi) \right| \lesssim \langle\xi\rangle ^{m-\rho|\alpha|+\delta|\beta|},
\end{equation*}
for all multi-indices $\alpha$ and $\beta$ and $(x,\xi)\in \Rl^n\times \Rl^n$, where 
$\langle\xi\rangle:= (1+|\xi|^2)^{1/2}.$\\

\noindent We refer to $m$ as the order of the amplitude and $\rho, \, \delta$ as its type. We will refer to the class $S_{\rho,\delta}^{m}(\Rl^n)$ with $0<\rho\leq 1$, $0\leq \delta<1$ as classical, to the class $S_{0,\delta}^{m}(\Rl^n)$ with $0\leq\delta<1$ as the exotic class, and to $S_{\rho,1}^{m}(\Rl^n)$ with $0\leq\rho\leq 1$ as the forbidden class of amplitudes.
\end{Def}

It turns out (as we will see later) that for the Besov-Lipschitz estimates we do not require any regularity in the $x$-variable and therefore we introduce the following class which was first defined by Kenig and Staubach in \cite{KS}.

\begin{Def}\label{symbol class Sm no regularity}
Let $m\in \Rl$ and $\rho\in [0,1]$. An \textit{amplitude} \emph{(}symbol\emph{)} $a(x,\xi)$ in the class $L^\infty S^m_{\rho}(\Rl^n)$ is a function $a:\Rl^n\times \Rl^n \to \Rl^n$ that verifies the estimate
\begin{equation*}
\Vert\partial_{\xi}^\alpha a(x,\xi) \Vert_{L^\infty_x(\Rn)} \lesssim \langle\xi\rangle ^{m-\rho|\alpha|},
\end{equation*}
for all multi-indices $\alpha$ and $\beta$ and $(x,\xi)\in \Rl^n\times \Rl^n$. Thus we only assume measurability in the $x$-variable.
\end{Def}
In our treatment of oscillatory integral operators, the phase functions take center stage, and the OIOs are classified according to their phases.
\begin{Def}\label{def:fk}
For $0<k<\infty$, we say that a real-valued \textit{phase function} $\varphi(x,\xi)$ belongs to the class $\textart F^k$, if
$\varphi(x,\xi)\in \mathcal{C}^{\infty}(\Rl^n \times\Rl^n \setminus\{0\})$ and satisfies the following estimates $($depending on the range of $k)$\emph:
\begin{itemize}
    \item for $k \geq 1$, 
        $$ 
|\partial^\alpha_\xi  (\varphi (x,\xi)-x\cdot\xi) |\leq
c_{\alpha} |\xi|^{k-1}, \quad  |\alpha| \geq 1 ,
$$
    \item for $0<k<1$, 
$$ 
|\partial^\alpha_\xi \partial^{\beta}_x  (\varphi (x,\xi)-x\cdot\xi) |\leq
c_{\alpha,\beta} |\xi|^{k-|\alpha|}, \quad |\alpha + \beta | \geq 1 ,
$$
\end{itemize}
for all $x\in \Rl^n$ and $|\xi|\geq 1$.

\end{Def}

 \begin{Rem}
     A well known and typical example of a phase in $\textart F^k$ is the phases $|\xi|^k+x\cdot\xi$ which are related to the operator $e^{i (\Delta)^{k/2}}$
 \end{Rem}
 An important condition on the phase function in the context of global regularity of OIOs, is the so-called strong non-degeneracy condition:
\begin{Def}\label{nondeg phase}
One says that the phase function $\varphi(x,\xi)$ satisfies the strong non-degeneracy condition \emph{(}or $\varphi$ is $\mathrm{SND}$ for short\emph{)} if
\begin{equation}\label{eq:SND}
	\big |\det (\partial^{2}_{x_{j}\xi_{k}}\varphi(x,\xi)) \big |
	\geq \delta,\qquad \mbox{for  some $\delta>0$ and all $(x,\xi)\in \mathbb{R}^{n} \times \Rl^n\setminus\{0\}$}.
\end{equation}
\end{Def}

In order the guarantee that our operators are globally $L^2$-bounded we should also put yet another condition of the phase which we shall henceforth simply refer to as the $L^2$-condition. The $L^2$-condition is essential for the validity of our main results. It allows us to control the $L^2$-behavior of the oscillatory integral operator.

\begin{Def}\label{L2 conditions}
One says that the phase function $\varphi(x,\xi)\in \mathcal{C}^{\infty}(\Rl^n \times\Rl^n)$ satisfies 
the weak $L^2$-condition if \begin{equation}\label{eq:L2 condition_new}
            |\partial_{x}^{\alpha} \partial_{\xi} \varphi(x,\xi )| \leq C_{\alpha}, \quad|\partial_{x} \partial_{\xi}^{\beta} \varphi(x, \xi)| \leq C_{\beta} 
            \end{equation}
            for $|\alpha|\geq 1$, $|\beta|\geq 1,$ all $x \in\mathbb{R}^{n} $ and $|\xi|\geq 1$.\\
We say that $\varphi(x, \xi)$ satisfies the strong $L^2$-condition if  \begin{equation}\label{eq:L2 condition_old}
            |\partial_{x}^{\alpha} \partial_{\xi}^{\beta} \varphi(x,\xi )| \leq C_{\alpha},
            \end{equation}
            for $|\alpha|\geq 1$, $|\beta|\geq 1,$ all $x \in\mathbb{R}^{n} $ and $|\xi|\geq 1$.
\end{Def}

\noindent Having the definitions of the amplitudes and the phase functions at hand, one has
\begin{Def}\label{def:OIO}
An oscillatory integral operator \emph{(OIO)} $T_a^\varphi$ with amplitude \linebreak$a\in S^{m}_{\rho, \delta}(\Rl^n)$ and a real valued phase function $\varphi$, is defined \emph{(}once again a-priori on $\mathscr{S}(\Rl^n)$\emph{)} by
\begin{equation}\label{eq:OIO}
T_a^\varphi f(x) := \int_{\Rl^n} e^{i\varphi(x,\xi)}\, a(x,\xi)\, \widehat f (\xi) \ddd\xi.
\end{equation}
If $\varphi\in \textart F^k$ and is $\mathrm{SND}$, then these operators will be referred to as oscillatory integral operators of order $k$.
\end{Def}

The $\mathrm{LF}(\mu)$-condition is a natural requirement which from the point of view of the applications into PDE's,  will always be satisfied and would not cause any loss of generality. It ensures that the operator behaves in a predictable and well-behaved manner, which is necessary for the analysis the low frequency portions of the operators.

\begin{Def}\label{Def:LFmu}
Assume that $\varphi(x,\xi)\in \mathcal{C}^{\infty}(\Rl^n \times \Rl^n \setminus \{0\}) $ is real-valued and $0<\mu\leq 1$. 
 We say that $\varphi$ satisfies the low frequency phase condition of order $\mu$, 
\linebreak\emph($\varphi$ satisfies \emph{LF}$(\mu)$-condition for short\emph),
if one has 
\begin{equation}\label{eq:LFmu}
|\partial^{\alpha}_{\xi}\partial_{x}^{\beta} (\varphi(x,\xi)-x\cdot \xi) |\leq c_{\alpha} |\xi|^{\mu-|\alpha|}, 
\end{equation}
for all $x\in \Rl^n$, $0<|\xi| \leq 2$ and all multi-indices 
$\alpha,\beta$.
\end{Def}

\begin{Rem}
As an example, note that the phase function associated to the water-wave equation is $x\cdot \xi+ |\xi|^{1/2}$, which satisfies the $\mathrm{LF}(\mu)$ condition with $\mu=\frac{1}{2}$. The phase function associated to the capillary wave equation is $x\cdot \xi+ |\xi|^{3/2},$  which is in $\mathrm{LF}(\mu)$ for any $\mu\in (0, 1]$.
\end{Rem}

We shall also need the following lemma to estimate the phase in the proofs for Sobolev-boundedness of OIOs with forbidden symbols (whenever the composition Theorem \ref{thm:left composition with pseudo} below is used), the proof of this lemma can be found in \cite[Lemma 4.1]{CISY}.

\begin{Lem}\label{lem:phaseestimates}
Assume $a(x,\xi)$ is an amplitude, and $\varphi(x,\xi)$ is a \emph{SND} phase function satisfying
\begin{equation}\label{lem:phaseestimatecond}
    |\d_\xi\d_x^\beta\varphi(x,\xi)|\leq c_\beta, \quad |\beta|\geq 1 \textnormal{ and } |\xi|\geq 1.
\end{equation}
Then for all $|\beta|\geq 1$, the following estimates
\begin{align}
    &|\nabla_x\varphi(x,\xi)|\sim |\xi|\label{lem:phaseestimate1}\\
    &|\d^{\beta}_x\varphi(x,\xi)|\lesssim \jap{\xi}\label{lem:phaseestimate2}
\end{align}
hold true for the phase function $\varphi$, on the support of $a(x, \xi)$, provided that the $\xi$-support of $a(x, \xi)$ lies outside the ball $B(0, R)$ for some large enough $R\gg1$ and $\d^{\beta}_x\varphi(x,\xi)\in L^\infty(\Rl^n\times \mathbb{S}^{n-1})$, for $|\beta|\geq 1.$
\end{Lem}
From this lemma it readily follows that
\begin{Cor}\label{bra fasfunktioner}
The phase functions in $\textart F^k$ that also satisfy the $L^2$-condition \eqref{eq:L2 condition_old}, all verify the estimates \eqref{lem:phaseestimate1} and \eqref{lem:phaseestimate2} of \emph{Lemma \ref{lem:phaseestimates}}.
\end{Cor}
\noindent We also recall a composition result proved in \cite{ATW}. This will be essential in the proof of the boundedness of OIOs with forbidden amplitudes on Sobolev spaces as well as in the proof of the Triebel-Lizorkin boundedness of oscillatory integral operators and in other situations.

\begin{Th}[\cite{ATW}]\label{thm:left composition with pseudo}
Let $m,s\in \Rl$, $\rho\in[0,1]$, $\delta\in[0,1)$. Suppose that \linebreak$ a(x, \xi)\in S_{\rho,\delta}^m (\Rl^n)$, $b(x,\xi)\in S^s_{1,0}(\Rl^n)$ and $\varphi(x,\xi)$ is a phase function that is smooth on $\supp a$ and verifies the conditions
\begin{align}\label{eq:composition_conditions}
    &|\xi| \lesssim |\nabla_x \varphi(x, \xi)| \lesssim |\xi|,  &|\partial_\xi^\alpha \partial _x^\beta \varphi (x, \xi)| \lesssim  \jap{\xi}^{1-|\alpha|}
\end{align}
for all $(x, \xi) \in \supp a$ and for all $|\alpha| \geq 0,$ and all $|\beta| \geq 1$.
\noindent For  $0<t\leq 1$ consider the parameter-dependent pseudodifferential operator 
\begin{equation*}
b(x, tD)f(x) := \int_{\Rl^n} e^{ix\cdot \xi}\,b(x,t\xi)\,\widehat f(\xi) \ddd \xi,
\end{equation*}
and the oscillatory integral operator
$$T_{a}^\varphi f(x) := \int_{\Rl^n} e^{i\varphi(x,\xi)}\,a(x, \xi)\,\widehat f(\xi) \ddd\xi.$$
Let $\sigma_t$ be the amplitude of the composition operator 
$T_{\sigma_t}^\varphi
:=b(x, tD)T_{a}^\varphi$ 
given by
\begin{equation*}
\sigma_t(x, \xi) := \iint_{\Rl^n\times \Rl^n} a(y, \xi)\, b (x,t\eta)\,e^{i(x-y)\cdot \eta+i\varphi(y,\xi)-i\varphi(x,\xi)} \ddd\eta \dd y.
\end{equation*}
Then for any $M\geq 1$ and all $0<\eps<1-\max(\delta,1/2)$, one can write $\sigma_t$ as

\begin{equation}\label{asymptotic expansion}
\sigma_t(x, \xi) =b(x,t\nabla_{x}\phase(x,\xi))\,a(x,\xi) + \sum_{0<|\alpha| < M}\frac{t^{\eps|\alpha|}}{\alpha!}\, \sigma_{\alpha}(t,x,\xi)+t^{\eps M} r(t,x,\xi),
\end{equation}

where, for all multi-indices $\beta, \gamma$ one has
\begin{align*}
|\partial^{\beta}_{\xi} \partial^{\gamma}_{x}\sigma_{\alpha}(t,x,\xi)| & \lesssim  t^{\min(s,0)}  \bra{\xi}^{s+m-(1-\max(\delta,1/2)-\varepsilon)|\alpha|-\rho|\beta|+\delta|\gamma|},  \\|\partial^{\beta}_{\xi} \partial^{\gamma}_{x} r(t,x,\xi)| &  \lesssim t^{\min(s,0)} \bra{\xi}^{s+m-(1-\max(\delta,1/2)- \varepsilon)M-\rho|\beta|+\delta|\gamma|}.
\end{align*}
\end{Th}

\begin{proof}
For a proof of this result see \cite{ATW}.
\end{proof}

\subsection{Some facts from the theory of classical function spaces}
\quad\\

\noindent 
Classical function spaces, such as Triebel-Lizorkin spaces and Besov-Lipschitz spaces, are central to the study of partial differential equations and the analysis of functions. In this section, we review some basic facts and results from the theory of classical function spaces that will be used throughout the rest of our paper.\\

We begin by introducing the definitions of classical Littlewood-Paley operators, Triebel-Lizorkin spaces, and Besov-Lipschitz spaces and discussing some of their key properties. Next, we present some important lemmas and theorems that provide useful estimates and bounds for functions in these spaces. We also introduce the Hardy-Littlewood maximal function, and how this operator can be leveraged in the analysis of various mathematical problems that will arise in the subsequent sections.

\begin{Def}\label{def:LP}
Let $\psi_0 \in \mathcal C_c^\infty(\Rl^n)$ be equal to $1$ on $B(0,1)$ and have its support in $B(0,2)$. Then let
$$\psi_j(\xi) := \psi_0 \left (2^{-j}\xi \right )-\psi_0 \left (2^{-(j-1)}\xi \right ),$$
where $j\geq 1$ is an integer and $\psi(\xi) := \psi_1(\xi)$. Then $\psi_j(\xi) = \psi\left (2^{-(j-1)}\xi \right )$ and one has the following Littlewood-Paley partition of unity

\begin{equation*}
    \sum_{j=0}^\infty \psi_j(\xi) = 1, \quad \text{\emph{for all }}\xi\in\Rl^n .
\end{equation*}

\noindent It is sometimes also useful to define a sequence of smooth and compactly supported functions $\Psi_j$ with $\Psi_j=1$ on the support of $\psi_j$ and $\Psi_j=0$ outside a slightly larger compact set. One could for instance set
\begin{equation*}
\Psi_j := \psi_{j+1}+\psi_j+\psi_{j-1},
\end{equation*}
with $\psi_{-1}:=\psi_{0}$.
\end{Def} \hspace*{1cm}\\
In what follows we define the Littlewood-Paley operators by

\begin{equation*}
    \psi_j(D)\, f(x)= \int_{\Rl^n}   \psi_j(\xi)\,\widehat{f}(\xi)\,e^{ix\cdot\xi} \ddd\xi,
\end{equation*}

where $\ddd \xi$ denotes the normalised Lebesgue measure ${\dd \xi}/{(2\pi)^n}$ and
\begin{equation*}
	\widehat{f}(\xi)=\int_{\Rl^n} e^{-i x\cdot\xi}\,f(x) \dd x,
\end{equation*}
is the Fourier transform of $f$.\\

Define $\|\cdot\|_{L^p(\ell^q)}$ and $\|\cdot\|_{\ell^q(L^p)}$ to be the quasi-norms
\begin{align*}
    \|(f_k)\|_{L^p(\ell^q)}&:=\Big\|\sum_{k} |f_k(\cdot)|^q\Big\|_{L^p(\Rn)}^{1/q},\\
    \|(f_k)\|_{\ell^q(L^p)}&:=\Big(\sum_{k} \|f_k\|_{L^p(\Rn)}^q\Big)^{1/q}.
\end{align*}
Using the Littlewood-Paley decomposition of Definition \ref{def:LP} we define the \emph{Triebel-Lizorkin space} $F^s_{p,q}(\Rl^n)$ and \emph{Besov-Lipschitz space} $B^s_{p,q}(\Rl^n)$.

\begin{Def}\label{def:TLspace}
	Let $s \in {\Rl}$ and $0< p <\infty$, $0< q \leq\infty$. The Triebel-Lizorkin space is defined by
	\[
	F^s_{p,q}(\Rl^n)
	:=
	\Big\{
	f \in {\mathscr{S}'}(\Rl^n) \,:\,
	\|f\|_{F^s_{p,q}(\Rl^n)}
	:=
	\|(2^{jqs}\psi_j(D) f)\|_{L^p(\ell^q)}<\infty
	\Big\},
	\]
where $\mathscr{S}'(\Rl^n)$ denotes the space of tempered distributions.	
\end{Def}
 
\begin{Def}\label{def:Besov}
	Let $0 < p,q \le \infty$ and $s \in {\mathbb R}$. The Besov-Lipschitz spaces are defined by
	\[
	{B}^s_{p,q}(\Rl^n)
	:=
	\Big\{
	f \in {{\SS}'(\Rl^n)} \,:\,
	\|f\|_{{B}^s_{p,q}(\Rl^n)}
	:=
	\|(2^{jqs}\psi_j(D) f)\|_{\ell^q(L^p)}<\infty
	\Big\}.
	\]
\end{Def}

\begin{Rem}\label{rem:TLspace}
Different choices of the basis $\{\psi_j\}_{j=0}^\infty$ give equivalent norms of $F^s_{p,q}(\Rl^n)$ in \emph{Definition \ref{def:TLspace},} see e.g. \cite{Triebel1}. We will use either $\{\psi_j\}_{j=0}^\infty$ or  $\{\Psi_j\}_{j=0}^\infty$ to define the norm of $F^s_{p,q}(\Rl^n)$.
\end{Rem}

\noindent We note that for $p=q=\infty$ and $0<s\leq 1$ we obtain the familiar Lipschitz space $\Lambda^s(\Rl^n)$, i.e. $B^s_{\infty,\infty}(\Rl^n)= \Lambda^s(\Rl^n)$. For $-\infty <s<\infty$ and $1\leq p<\infty,$ $F^s_{p,2}(\Rl^n)=H^{s,p}(\Rl^n)$ (various $L^p$-based Sobolev and Sobolev-Slobodeckij spaces) and for $0<p<\infty$, $F^0_{p,2}(\Rl^n)=h^p(\Rl^n)$ (the local Hardy spaces). Moreover the dual space of $F^{0}_{1,2}(\Rl^n)$ is $\mathrm{bmo}$ (the local version of $\mathrm{BMO}$).\\

The following estimate will be useful in the $L^1(\Rn)$ to $bmo(\Rn)$ boundedness of Oscillatory integral operators. Let $X^s_{p,q}(\Rn)$ be either $B^s_{p,q}(\Rn)$ or $F^s_{p,q}(\Rn)$. Then for all $-\infty <s<\infty$ and $0<p,q\leq \infty$ one has 
\begin{equation}\label{pointwise multiplier}
    \|fg\|_{X^s_{p,q}(\Rn)}\lesssim \Big(\sum_{|\alpha|\leq M}\sup_{x\in\Rn}|\d^\alpha f(x)|\Big)\|g\|_{X^s_{p,q}(\Rn)}.
\end{equation}
Two other useful facts which will be useful to us is that
for $-\infty <s<\infty$ and $0<p\leq \infty$ one has
    \begin{equation}\label{equality of TL and BL}
        B^s_{p,p}(\Rl^n)= F^s_{p,p}(\Rl^n),
    \end{equation}
    and
    \begin{equation}\label{embedding of TL}
        F^{s+\varepsilon}_{p,q_0}(\Rl^n)\xhookrightarrow{} F^s_{p ,q_1}(\Rl^n)\quad\text{and}\quad
        B^{s+\varepsilon}_{p,q_0}(\Rl^n)\xhookrightarrow{} B^s_{p ,q_1}(\Rl^n) 
    \end{equation}
for $-\infty <s<\infty$, $0<p< \infty$, $0<q_0,q_1 \leq \infty$ and all $\varepsilon>0$.\\

Furthermore, for $s'\in \Rl$, the operator $ (1-\Delta)^{s'/2}$ maps ${F}^s_{p,q}(\Rl^n)$ isomorphically into ${F}^{s-s'}_{p,q}(\Rl^n)$ and ${B}^s_{p,q}(\Rl^n)$ isomorphically into ${B}^{s-s'}_{p,q}(\Rl^n),$ see \cite[p. 58]{Triebel1}.\\

In connection to estimates for linear operators in Triebel-Lizorkin spaces, one often encounters  the well-known Hardy-Littlewood's maximal function
 \begin{equation*}\label{HLmax}
  \mathcal{M} f(x):=\sup _{B\ni x}\frac{1}{|B|} \int_{B}|f(y)| \dd y,   
 \end{equation*} where the supremum is taken over all balls $B$ containing $x$. For \(0<p<\infty\), one also defines 
$$(\mathcal{M}_{p} f(x):=\left(\mathcal{M}\left(|f|^{p}\right)\right)^{1 / p}.$$

We now present the following abstract lemma which we will make use of in relation to the boundedness of OIOs with amplitudes in the class $S^m_{1,1}(\Rn)$.
\begin{Lem}\label{Usingequivalenceofnormslemma}
    Let $X^{s}_{p,q}$ be either $F^{s}_{p,q}(\Rn)$ or $B^{s}_{p,q}(\Rn)$. Assume that $\{u_k\}_{k\geq 0}\in X^{s}_{p,q}$ is an arbitrary sequence such that $\big\| \{2^{ks}u_k\}_{k=0}^\infty\big\|_{_{L^p(l^q)}}\lesssim \Vert f\Vert_{X^{s}_{p,q}}$. Moreover let $\{h_{k,l}(x)\}_{k,l\geq 0}$ be a sequence such that the spectrum of $h_{k,l}(x)$ is in $B(0,2^{k+l+2})$ and that there is some sufficiently large $N$ such that $h_{k,l}$ satisfies the pointwise estimate,
    \begin{equation}\label{piecewiseest}
        |h_{k,l}(x)|\lesssim 2^{-Nl}\,\mathcal{M}_{r}u_k(x),\quad\exists r \in (0,\min\{p,q\}).
    \end{equation}
    Let $g_l=\sum_{k\geq 0}h_{k,l}$.
    If $0<p,q<\infty$ and $s>n(\frac{1}{\mathrm{min}( p, q,1)}-1)$ then
    \begin{equation}
        \norm{\sum_{l\geq 0} g_l }_{X^{s}_{p,q}} \lesssim \Vert f\Vert_{X^{s}_{p,q}}.
    \end{equation}
\end{Lem}

\begin{proof}
For a proof see \cite{ATW}.
\end{proof}

Another important and useful fact about Besov-Lipschitz and Triebel-Lizorkin spaces is the following:

\begin{Th}\label{thm:invariance thm}
Let $\eta: \Rn \to\Rn$ with $\eta(x)=(\eta_1 (x), \dots,\eta_n (x))$ be a diffeomorphism such that $\abs{\det D\eta (x)}\geq c>0$, $\forall x\in \Rn$ \emph{(}$D\eta$ denotes the Jacobian matrix of $\eta$\emph{)}, and $\Vert\partial^{\alpha}\eta_j (x)\Vert_{L^\infty(\Rn)}\lesssim 1$ for all $j\in \{1,\dots, n\}$ and $|\alpha|\geq 1.$ Then
for $s\in \Rl$, $0<p<\infty$ and $0<q\leq \infty$ one has $$\Vert f\circ \eta\Vert_{F^{s}_{p,q}(\Rn)}\lesssim \Vert f\Vert_{F^{s}_{p,q}(\Rn)}.$$

The same invariance estimate is also true for Besov-Lipschitz spaces $B^{s}_{p,q}(\Rn)$ for  $s\in \Rl$, $0<p\leq\infty$ and $0<q\leq \infty$.
\end{Th}
For a proof see J. Johnsen, S. Munch Hansen and W. Sickel \cite[Corollary 25]{JMHS}, and  H. Triebel \cite[Theorem 4.3.2]{Triebel2}. \\

The rest of this section is dedicated to setting up some definitions which will be used in relation to transference to Triebel-Lizorkin spaces (section \ref{subsec:OIO_lowP_transference}) and $h^p\to L^p$ estimates of oscillatory integral operators (section \ref{sec:hplp}).\\

For Definition \ref{Def:Qbe} and Definition \ref{def:hpatom}, let $[x]$ denote the integer part of $x$.
\begin{Def}\label{Def:Qbe}
For a closed cube $Q$, define
\begin{enumerate}
    \item[$(i)$] $c_Q$ as the centre of $Q$,
    \item[$(ii)$] $l_Q$ as the side length of $Q$,
    \item[$(iii)$]     $k_Q:=[1-\log_2(l_Q)]$, 
    \item[$(iv)$] $\chi_Q$ as the characteristic function of $Q$, i.e. $\chi_Q(x):=\diffcases{1,& x\in Q\\0,& x\notin Q}$.
    \item[$(v)$] $cQ$ as the cube with centre $c_Q$ and length $l_{c\,Q}=c\,l_Q$, where $c\in \Rl_{>0}$.
\end{enumerate}

\end{Def}
Observe that $k_Q$ is the unique integer such that $$2^{-k_Q}< l_Q\leq 2^{-(k_Q-1)}.$$
\begin{Def}\label{def:hpatom}
A function $\at$ is called a $h^p$-atom if there exists a cube $Q$ such that the following three conditions are satisfied:
\begin{enumerate}
\item[$(i)$] $\supp \at\subset Q$,
\item[$(ii)$] $\displaystyle \sup_{x\in Q}|\at(x)|\leq 2^{k_Qn/p},$
\item[$(iii)$] If $k_Q\geq 1$, $\displaystyle \mathfrak M_{\at}= \left[ n\brkt{\frac 1p -1} \right ]$, 
then $\displaystyle \int_{\Rn} x^{\alpha}\at(x)\dd x=0,$ for $|\alpha|\leq \mathfrak M_{\at}$. No further condition is assumed if $k_Q\leq 0.$ 
\end{enumerate}

It is well known $($see \cite{Triebel1}$)$ that a distribution $f\in h^p (\Rn)$ has an atomic decomposition
\nm{eq:hpsumdefinition}{
    f=\sum_{j=0}^\infty\lambda_{j}\at_{j},
}
where the $\lambda_{j}$ are constants such that $$\displaystyle \inf_{\set{\lambda_j}}\sum_{j=0}^\infty|\lambda_{j}|^{p}\sim\Vert f\Vert_{h^p(\mathbb{R}^{n})}^{p}=\Vert f\Vert_{F_{p,2}^0(\mathbb{R}^{n})}^{p}$$ and the $\at_{j}$ are $h^p$-atoms.\\

\end{Def}

In the analysis of the boundedness of OIOs on classical function spaces one typically decomposes the operator into a low- and a high-frequency part, where the low frequency part corresponds to an amplitude that is smooth and compactly supported in the $\xi$ variable. The following theorem, which was proven in \cite{CISY}, addresses the boundedness of the low-frequency portion of the OIOs.
\noindent The main boundedness result for the low frequency part of OIOs.

\begin{Th}[\cite{CISY}]\label{thm:low_freq_TL_BL_OIO}
For $\rho,\delta\in [0,1]$ let $a(x,\xi)\in S^m_{\rho,\delta}(\Rn)$ be an amplitude which is compactly supported in $\xi$. Suppose also that $0<q_1,q_2 \leq\infty$, $s_1,s_2\in \Rl$. If $\varphi (x,\xi)$ verifies the \emph{LF($\mu$)}-condition, then $T_{a}^\varphi $ is $F_{p,q_1}^{s_1}(\Rl^n)$ to $ F_{p,q_2}^{s_2}(\Rl^n)$ bounded, for $\frac{n}{n+\mu}<p\leq\infty$
Moreover, all the Triebel-Lizorkin estimates above may be replaced by the corresponding Besov-Lipschitz estimate.
\end{Th}

\begin{proof}
    See \cite[Lemma 6.3]{CISY}.
\end{proof}

\noindent Here we recall a theorem that allows one to lift $h^p\to L^p$ boundedness to Besov-Lipschitz boundedness. This result can be applied to a wide class of oscillatory integral operators. Observe that this lift is only valid for the classical amplitudes when $\delta<1$.

\begin{Th}[\cite{IRS}]\label{theorem:BL-booster theorem}
Let $0<p,q<\infty$, $m \leq 0$, $a\in S^m_{\rho,\delta}(\Rn)$ such that $a$ is supported in $\Rn\times\Rn\setminus B(0,1)$. Suppose $\varphi$ is a phase function that verifies the conditions \eqref{eq:composition_conditions} of \emph{Theorem \ref{thm:left composition with pseudo}}.
If $T^\varphi_a$ is $h^p\to L^p$-bounded, then $T^\varphi_a$ is bounded from $B^s_{p,q}\to B^s_{p,q}$ for all $s\in \Rl$. 
\end{Th}

\begin{proof}
    For a proof see \cite{IRS}.
\end{proof}

The following theorem due to M. Pramanik, K. Rogers and A. Seeger could be found in \cite[Theorem 2.1]{PRS}.
\begin{Th}[\cite{PRS}]\label{thm:PRS_proto}
Let $ 1<q < p < \infty$, and $ 0<b<n$. Assume that the sequence of operators $S_j$ satisfy
\begin{align}
&\sup_{
j>0}
2^{jb/p}\| S_j\|_{L^p\to L^p} \leq A_0,\label{eq:PRS1_proto}
\\&
\sup_{
j>0}
2^{jb/q}\|S_j\|_{L^q\to L^q} \leq B_0.\label{eq:PRS2_proto}
\end{align}
Furthermore, assume that for each cube $Q$ there is a measurable set $\mathcal E_{Q}$ and a constant $\Gamma\geq 1$ such that
\begin{equation}\label{ mattet}
|\mathcal E_{Q}| \leq \Gamma  \max\{l_{Q}^{n-b},l_{ Q}^n\},    
\end{equation}
where $l_{Q}$ is the side length of $ Q$ as in \emph{Definition} \emph{\ref{Def:Qbe}}, and assume further that for every $j\in \mathbb{N}$ and every cube \(Q\) with \(2^{j} l_{Q} \geq 1\), one has
\nm{eq:PRS}{
\sup_{x\in Q}
\int_{\Rn \setminus \mathcal E_{ Q}} 
|K_j(x, y)| \dd y \leq B_1  \max \{2^{-j\eps}l_{ Q}^{-\eps}, 2^{-j\eps}\},}
for some $\eps>0$.\\
Then if \begin{equation}
{B}_2:=B_{0}^{q / p}\left(A_0 \Gamma^{1 / p}+B_{1}\right)^{1-q / p}
\end{equation} one has
\eq{
\norm{ \brkt{\sum_{j=0}^\infty 
2^{jbr/p}|\Psi_j(D) S_jf_j|^r}^{1/r}}_{L^p(\Rn)} \lesssim A_0\left[\log \left(3+\frac{{B}_2}{A_0}\right)\right]^{1 / r-1 / p}  \brkt{ \sum_{j=0}^\infty \|f_j\|_{L^p(\Rn)}^p}^{1/p}.
}
where $\Psi\in \mathscr{S}(\Rl^n)$, $\Psi_j(D):=\Psi(2^{-j}D)$ and $f_j$ is a sequence of functions.
\end{Th}

\section{Basic kernel estimates for oscillatory integral operators}\label{basic kernel estimates}

\noindent In order to show the Besov-Lipschitz boundedness of OIOs related to forbidden amplitudes we prove some preliminary estimates on the kernel of the operators. What follows is a decomposition which, in contrast to the second angular-radii frequency decomposition of C. Fefferman  (see \cite{Feff}), decomposes the annuli of the standard Littlewood-Paley decomposition into balls. Furthermore, this decomposition generalizes the decomposition in \cite{CISY}, which used balls of constant radii to decompose the case corresponding to the OIOs of Schr\"odinger type phase functions in $\textart F^2$, to the more general case of balls of varying radii adapted to all OIOs with phases in the class $\textart F^k$.

\begin{figure}
\begin{tikzpicture}

        \draw (0,0) circle (1);
        \draw (0,0) circle (2);
        \filldraw[red] (0.6,1.04) circle (0.05);
    
        \draw (0.2,0.4) node[anchor=east] {$\xi_j^\nu$};
    
        \foreach \i in {20,40,...,360}
        {
            \draw[rotate = \i] (1.2,0) circle (0.3);
        }
        \foreach \i in {18,36,...,360}
        {
            \draw[rotate = \i] (1.5,0) circle (0.3);
        }
        \foreach \i in {15,30,...,360}
        {
            \draw[rotate = \i] (1.8,0) circle (0.3);
        }

        \draw[->] (0.15,0.6) -- (0.6,1.04);
\end{tikzpicture}
\caption{Covering of $\supp \psi_j\subset \{ 2^{j-1}\leq|\xi|\leq 2^{j+1}\}$ with balls of radii $2^{j(1-k/2)}$ and centres $\xi_j^\nu$. }\label{fig:schrodingercover}
\end{figure}
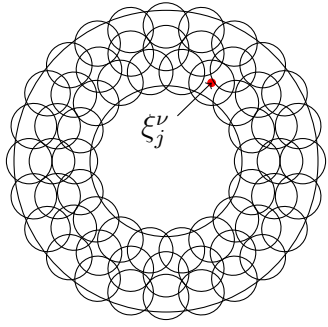

\begin{Def}\label{second decomposition} Let $k>0$.
    We make the following decomposition of the integral kernel 
$$K(x,y)= \int_{\Rl^n}  a(x,\xi)\,e^{i\varphi(x,\xi)-iy\cdot \xi}\,  \ddd \xi.$$
We introduce a standard Littlewood-Paley partition of unity $\sum_{j=0}^\infty \psi_j(\xi) =1$ with $\supp \psi_0 \subset B(0,2)$, and $\supp \psi_j\subset \{ 2^{j-1}\leq|\xi|\leq 2^{j+1}\}$ for $j\geq 1$ as in \emph{Definition \ref{def:LP}}. \\

Then for every $j\geq 0$ we cover  
$\supp \psi_j$ with open balls $C_j^\nu$ with radii $2^{j(1-k)}$ and centres $\xi_j^\nu$, where $\nu$ runs from $1$ to $\mathscr{N}_j:=O(2^{njk})$. See \emph{Figure \ref{fig:schrodingercover}} for an illustration. Observe that $|C_j^\nu|\lesssim 2^{n j(1-k)}$ uniformly in $\nu$. 
Now take $u\in \mathcal{C}_c^{\infty}(\Rl^n)$, with $0\leq u\leq 1$ and supported in $B(0,2)$ with $u=1$ on $\overline{B(0,1)}.$\\

Next set
$$\chi_j^\nu(\xi) := \frac{u(2^{-(1-k) j}(\xi -\xi_{j}^{\nu}))} {\sum_{\kappa=1}^{\mathscr{N}_j} u(2^{-(1-k) j}(\xi-\xi_{j}^{\kappa}))},$$
and note that
$$\sum_{j=0}^\infty \sum_{\nu=1}^{\mathscr N_j} \chi_j^\nu(\xi)\,\psi_j(\xi) =1.$$
Now we define the second frequency localized pieces of the kernel above as
\begin{equation*}
K_{j}^\nu (x,y) := \int_{\Rl^n}  \psi_j(\xi)\,\chi_j^\nu(\xi)\,e^{i\varphi(x,\xi)-iy\cdot \xi}\,a(x,\xi)  \ddd \xi.
\end{equation*}
\end{Def}

\begin{Lem}\label{kernel lemma}
       Let $0\leq \rho\leq 1$, $0\leq \delta< 1$, $n\geq 1$, and $0<k<\infty $. Assume that $\varphi\in \textart F^k$ {is \emph{SND}}, satisfies the $L^2$-condition \eqref{eq:L2 condition_old}. 
Then if $a(x,\xi)\in L^{\infty}S^{m}_{\rho}(\Rl^n)$ $($see \emph{Definition \ref{symbol class Sm no regularity}}$)$, and $K_{j}^\nu (x,y)$ is defined as in \emph{Definition \ref{second decomposition}}, then we have
\begin{equation}\label{main kernel estimate for AAW 1}
|\partial_y^\alpha K_j^\nu (x,y)|
\lesssim  \frac {2^{j(m+|\alpha|)} 2^{jn(1-k)}}{\jap{2^{jw(k,\rho)}(\nabla_\xi \varphi(x,\xi_j^\nu)-y)}^{M}},
\end{equation}
for all multi-indices $\alpha$, all $j\geq 0$ and $M\geq 0$ and where
$$w(k,\rho)=
\begin{cases}
    \min\{\rho,1-k\} & 0<k<1,\\
    k-1 & k\geq1.\\
\end{cases}$$
\end{Lem}

\begin{proof}

Observe that we have
\begin{equation*}
    |\partial^{\gamma}_\xi \chi_j^{\nu}(\xi)|\lesssim 2^{j(k -1)|\gamma|}
\end{equation*}

Therefore, for any multi-index $\alpha$ and any $j\geq 0$ we have
\begin{align*}
\partial^{\alpha}_{y}K_{j}^{\nu}(x,y) 
& = \int_{\Rl^n}  \psi_j(\xi)\,\chi_j^\nu(\xi)\,
\partial^{\alpha}_{y} e^{i(\varphi(x,\xi )-y\cdot \xi)} \,a(x,\xi) \ddd \xi \\
& = 
\int_{\Rl^n}
e^{i(\varphi(x,\xi )-y\cdot \xi)}\, \sigma_j^{\alpha,\nu}(x,\xi) \ddd\xi,
\end{align*}

where
$$\sigma_{j}^{\alpha,\nu}(x,\xi)
:=  \psi_j(\xi)\,\chi_j^\nu(\xi)\,(-i\xi)^\alpha \, a(x,\xi).
$$
Using the assumption that $a(x,\xi)\in L^{\infty}S^{m}_{\rho}(\Rl^n)$, we deduce that for any multi-index $\gamma$, any $j\geq 0$ and any $\nu$ one has
\begin{equation}\label{amplitude derivative estim}
 |\partial^{\gamma}_\xi \sigma_j^{\alpha,\nu}(x,\xi)|\lesssim 2^{j(m+|\alpha|-\min\{\rho,1-k\} |\gamma|)},
\end{equation}

If we now set $\vartheta_j^\nu(x,\xi):=\varphi(x,\xi)-\xi \cdot \nabla_\xi \varphi(x,\xi_j^\nu)$, then we can write
\begin{equation*}
\partial^{\alpha}_{y}K_{j}^\nu(x,y) = \int_{\Rl^n}  e^{i(\nabla_\xi\varphi(x,\xi_j^\nu )-y)\cdot \xi}\,e^{i\vartheta_j^\nu(x,\xi)}\,\sigma_{j}^{\alpha,\nu}(x,\xi)  \ddd \xi.
\end{equation*}
 Now we estimate the derivatives of $\vartheta$ in $\xi$ on the support of $\sigma_{j}^{\alpha,\nu}(x,\xi)$. To this end, the mean-value theorem and Definition \ref{def:fk} yields for $k\geq1$
\eq{
|\partial_{\xi_l}\vartheta_j^\nu(x,\xi)|
&=|\partial_{\xi_l}\varphi(x,\xi)- \partial_{\xi_l} \varphi(x,\xi_j^\nu)|   
=  \Big| (\xi-\xi_j^\nu)\cdot \int_0^1(\nabla_\xi \partial_{\xi_l}\varphi)(x,t\xi+(1-t)\xi_j^\nu)\dd t \Big| \\
&\lesssim 2^{j(1-k)} \sup_{t\in [0,1]} |t\xi+(1-t)\xi_j^\nu)|^{k-1} \\
&\lesssim 2^{j(1-k)}2^{j(k-1)}= 1
}
and
\eq{
|\partial_{\xi}^{\alpha}\vartheta_j^\nu(x,\xi)| = |\partial_{\xi}^{\alpha}\varphi(x,\xi)|\lesssim 2^{j|\alpha|\max\{(k-1),0\}},\quad \text{ for all } |\alpha| \geq 2.
}
Therefore by Faa di Bruno's formula we obtain that
\begin{align}\label{linearizeddecayestimate}
    |\partial_{\xi}^{\gamma} e^{i  \vartheta^{\nu}_{j}(\xi,x)}| \lesssim
    \sum_{\gamma_1+ \dots+ \gamma_r=\gamma} |\partial_{\xi}^{\gamma_1} \vartheta^{\nu}_j| \dots |\partial_{\xi}^{\gamma_r} \vartheta^{\nu}_j|
    \lesssim 2^{j\max\{(k-1),0\}|\gamma|}.
\end{align}

Observe the simple estimate
\begin{equation}
   |\partial^{\alpha}_{y}K_{j}^\nu(x,y) |\lesssim 2^{j(m+|\alpha|)} 2^{jn(1-k)}
\end{equation}
together with integration by parts yield
\begin{equation}
   |\partial^{\alpha}_{y}K_{j}^\nu(x,y) |
\lesssim \frac {2^{j(m+|\alpha| +w(k,\rho) M)} 2^{jn(1-k)}}{|\nabla_\xi \varphi(x,\xi_j^\nu)-y|^{M}}.
\end{equation}
where $w(k,\rho)=\max\{-\min\{\rho,1-k\},\max\{(k-1),0\}\}$. These equalities, the observation that $\supp \sigma_{j}^{\alpha,\nu}\subset C^{\nu}_{j},$ with $|C^{\nu}_j|=O(2^{jn(1-k)})$ uniformly in $\nu$ and $j$, the estimates for the derivatives of $\vartheta$, and \eqref{amplitude derivative estim} yield 
\begin{equation*}
|\partial_y^\alpha K_j^\nu (x,y)|
\lesssim \frac {2^{j(m+|\alpha| )} 2^{jn(1-k)}}{\jap{2^{jw(k,\rho)}(\nabla_\xi \varphi(x,\xi_j^\nu)-y)}^{M}}\lesssim \frac {2^{j(m+|\alpha| )} 2^{jn(1-k)}}{\jap{2^{jw(k,\rho)}(\nabla_\xi \varphi(x,\xi_j^\nu)-y)}^{M}},
\end{equation*}
for all multi-indices $\alpha$, all $j\geq 0$ and $M\geq 0$.
\end{proof}

\section{\texorpdfstring{$L^2$}{}-boundedness of Oscillatory integral operators}\label{L2 section}
Traditionally, the way to prove $L^2$-boundedness results for oscillatory integral operators at the endpoint of the range of exponents has been to utilize the almost orthogonality principle, as established in the continuous version of the Cotlar-Stein lemma by Calderón and Vaillancourt, see i.e. \cite{CV} for a proof. This lemma, which was originally developed by Cotlar and Stein, is a powerful tool for establishing bounds on oscillatory integral operators.

\begin{Lem}\label{calderon-vaillancourt lemma}
Let $\mathscr{H}$ be a Hilbert space, and $A(\xi)$ a family of bounded linear endomorphisms of $\mathscr{H}$ depending on $\xi \in \Rn .$
Assume the following three conditions hold:
\begin{enumerate}
  \item[$(i)$] the operator norm of $A(\xi)$ is less than a number $C$ independent of $\xi.$
  \item[$(ii)$] for every $u\in \mathscr{H}$ the function $\xi\mapsto A(\xi)u$ from $\Rn \mapsto \mathscr{H}$ is continuous for the norm topology of $\mathscr{H}.$
  \item[$(iii)$] for all $\xi_{1}$ and $\xi_2$ in $\Rn $
  \begin{equation}\label{cotlar estimates}
    \Vert A^{\ast}(\xi_1) A(\xi_2)\Vert \leq h(\xi_1 , \xi_2)^2,\,\,\,\text{and}\,\,\,  \Vert A(\xi_1) A^{\ast}(\xi_2)\Vert \leq h(\xi_1 , \xi_2)^2 ,
  \end{equation}
  with $h(\xi_1 , \xi_2 )\geq 0$ is the kernel of a bounded linear operator on $L^2$ with norm $K$.
\end{enumerate}
Then for every $E\subset \Rn $, with $|E|<\infty$, the operator $A_{E}=\int_{E} A(\xi)  \ddd\xi$ defined by $\langle A_{E} u, v\rangle _{\mathscr{H}} = \int_{E} \langle A(\xi) u, v\rangle _{\mathscr{H}} \ddd\xi,$ is a bounded linear operator on $\mathscr{H}$ with norm less than or equal to $K.$
\end{Lem}
Another useful fact that will aid us in the estimate of the oscillatory integrals is the following lemma whose proof could be found in \cite{DS}.\\
\begin{Lem}\label{integration by parts lem}
Let $s(x)$ and $F(x)$ be real-valued smooth functions in $\Rl^n$, and
\begin{equation}\label{definition of integration by parts}
  Lu(x):= D^{-2} (1-i s(x)\langle\nabla_{x}F, \nabla_{x}\rangle)u(x),
  \end{equation}
    with $D:= (1+ s(x)|\nabla_{x} F|^{2})^{1/2}.$ Then
  \begin{enumerate}
    \item[$(i)$] $L (e^{iF(x)})= e^{iF(x)}$
    \item[$(ii)$] if ${}^{t}L$ denotes the formal transpose of $L,$ then for any positive integer $N,$ $({}^{t}L)^{N} u(x)$ is a finite linear combination of terms of the form
    \begin{equation}\label{mainterm}
      CD^{-k} \Big\{\prod_{\mu=1}^{p}\partial^{\alpha_\mu}_{x} s(x)\Big\}\Big\{\prod_{\nu =1}^{q}\partial^{\beta_{\nu}}_{x} F(x)\Big\}\partial^{\gamma}_{x} u(x),
    \end{equation}
    with
    \begin{multline}\label{relations for order of derivatives}
      2N\leq k\leq 4N ;\, k-2N\leq p\leq k-N ;\,
      |\alpha_{\mu}|\geq 0;\, \sum_{\mu=1}^{p} |\alpha_{\mu}|\leq N\\
       k-2N\leq q\leq k-N;\, |\beta_{\nu}|\geq 1;\,\sum_{\mu=1}^{q} |\beta_{\nu}|\leq q+N;\, |\gamma| \leq N.
    \end{multline}
  \end{enumerate}
\end{Lem}

\begin{Th}\label{thm:Calderon-Vaillancourt for OIOs}
If $m=\min(0,\frac{n}{2}(\rho-\delta)),$ $0\leq \rho\leq 1$, $0\leq \delta<1,$ $a\in S^{m}_{\rho, \delta}(\Rn)$ and $\varphi\in \mathcal{C}^{\infty}(\mathbb{R}^n \times \mathbb{R}^n)$ is \emph{SND} for all $(x, \xi)\in \mathbb{R}^n \times \mathbb{R}^n$ and satisfies the weak $L^2$-condition \eqref{eq:L2 condition_new} on the support of $a(x,\xi)$, then the operator $$T^{\varphi}_{a} f(x)= \int_{\Rn} a(x,\xi)\, e^{i\varphi(x,\xi)} \hat{f}(\xi) \ddd\xi$$ is bounded on $L^2(\mathbb{R}^n).$
\end{Th}

\begin{proof}
Let $\chi(x, \xi)\in \mathcal{C}_{c}^\infty(\Rl ^n \times \Rl^n)$ be such that $\chi(x, \xi)=1$ for $|x|^2 + |\xi|^2\leq 1$ and define $a_\varepsilon(x, \xi):= \chi(\varepsilon x, \varepsilon \xi) \, a(x, \xi).$ Since $a_\varepsilon$ converges to $a$
in $S^{m}_{\rho,\delta}(\Rl^n)$ so that for any $f\in \mathscr{S}(\Rl ^n)$ ,  $T_{a_\varepsilon}^{\varphi} f$ converges, as $\varepsilon\to 0$, to $T_{a}^{\varphi} f$ in $\mathscr{S}(\Rl^n)$. Since the seminorms of $a_\varepsilon$ are bounded by a constant (depending on $\chi$) times the seminorms of $a$, we can therefore assume from now on that $a(x, \xi)\in \mathcal{C}_{c}^\infty(\Rl ^n \times \Rl^n ).$ Later on, of course, the estimates that we obtain won't depend on the support of $a$.\\

Furthermore, we observe that since for $\delta \leq \rho$, $S^{0}_{\rho, \delta}(\Rn) \subset S^{0}_{\rho, \rho}(\Rn),$ it is enough to show the theorem for $0\leq \rho\leq \delta<1$ and $m=\frac{n}{2}(\rho-\delta).$ Using the unitarity of the Fourier transform in $L^2(\mathbb{R}^n)$ and a $TT^{\ast}$ argument, it is enough to show that the operator
\begin{equation}\label{Tb amplitude presentation}
  T_{b} f(x)= \iint b(x,y,\xi)\, e^{i\varphi(x,\xi)-i\varphi(y,\xi)}\, f(y) \dd y  \ddd\xi,
\end{equation}
where $b(x, y, \xi):=a(x, \xi)\overline{a(y, \xi)}$ satisfies the estimate
\begin{equation}\label{derivtives of b}
|\partial^{\alpha}_{\xi}\partial^{\beta}_{x} \partial^{\gamma}_{y} b(x,y,\xi)|\leq C_{\alpha\, \beta\,\gamma} \langle \xi \rangle ^{m_1-\rho|\alpha|+\delta(|\beta| +|\gamma|)},
\end{equation}
 with $m_1=n(\rho-\delta)$ and $0\leq \rho\leq \delta<1,$ is bounded on $L^2(\mathbb{R}^n).$ Moreover due to assumption of compact support of $a(x, \xi)$ we can also assume that $b\in \mathcal{C}_{c}^\infty(\Rl ^n \times \Rl^n \times \Rl^n)$, under the understanding that the norm estimates that we obtain will be independent of the support of $b$.\\
Now we introduce a differential operator
   $$ L:= D^{-2} \Big\{1-i\langle \xi\rangle^{\rho}\big(\langle\nabla_{\xi}\varphi(x,\xi)-\nabla_{\xi}\varphi(y,\xi), \nabla_{\xi}\rangle\big)\Big\},$$
with $D=(1+\langle \xi\rangle ^{\rho}|\nabla_{\xi}\varphi(x,\xi)-\nabla_{\xi}\varphi(y,\xi)|^2)^{\frac{1}{2}}.$ It follows from Lemma \ref{integration by parts lem} that
 \begin{equation*}
  L (e^{i\varphi(x,\xi)-i\varphi(y,\xi)})= e^{i\varphi(x,\xi)-i\varphi(y,\xi)}
\end{equation*}
and that $({}^t L)^{N}b(x,y,\xi)$ is a finite sum of terms of the form
\begin{equation}\label{eq:finite_sum_of_terms}
D^{-k} \bigg\{\prod_{\mu=1}^{p}\partial^{\alpha_\mu}_{\xi} \langle \xi\rangle^{\rho}\bigg\}\bigg\{\prod_{\nu =1}^{q} \big(\partial^{\beta_{\nu}}_{\xi}\varphi(x,\xi)-\partial^{\beta_{\nu}}_{\xi}\varphi(y,\xi)\big)\bigg\}\,\partial^{\gamma}_{\xi}b(x,y,\xi),
\end{equation}
with $p$, $\alpha_\mu$, $\beta_\nu$, $q$ and $\gamma$ quantified by \eqref{relations for order of derivatives}.
Furthermore since $\varphi$ is SND, we can use Proposition 1.11 in \cite{DS} to show that
\begin{equation}
\label{lowerbound for phi in x}
|\nabla_{\xi} \varphi(x,\xi)-\nabla_{\xi}\varphi(y,\xi)|\geq c_1 |x-y|
\end{equation}
 \begin{equation}
\label{lowerbound for phi in xi}
|\nabla_{z} \varphi(z,\xi_1)-\nabla_{z}\varphi(z,\xi_2)|\geq c_2|\xi_1-\xi_2|.
\end{equation}
Using \eqref{lowerbound for phi in x}, \eqref{relations for order of derivatives}, \eqref{derivtives of b} and \eqref{eq:finite_sum_of_terms}, we have 
\begin{equation}\label{x derivative estimate}
  |\partial^{\sigma}_{x}({}^t L)^{N}b(x,y,\xi)|\leq C \Lambda(\langle \xi\rangle ^{\rho} (x-y))\langle \xi\rangle ^{m_1+\delta|\sigma|},
\end{equation}
where $\Lambda$ is a  function with $\int_{\Rn} \Lambda(\xi)  \ddd \xi\lesssim 1$.
 Integration by parts using $L$, $N$ times, in \eqref{Tb amplitude presentation} one has
 \begin{equation}
 T_{b} f(x)= \iint_{\Rn\times\Rn} c(x,y,\xi)\, e^{i\varphi(x,\xi)-i\varphi(y,\xi)}\, f(y) \dd y  \ddd\xi,
 \end{equation}
 with $c(x,y,\xi)=({}^t L)^{N}b(x,y,\xi)$ and
 \begin{equation}\label{derivative estimates for c}
 |\partial^{\sigma}_{x}c(x,y,\xi)|\leq C \Lambda(\langle \xi\rangle ^{\rho} (x-y))\langle \xi\rangle ^{m_1+\delta|\sigma|}
 \end{equation}
 and the same estimate is valid for $\partial^{\sigma}_{y}c(x,y,\xi).$
 From this we get the representation
 \begin{equation}
   T_{b} = \int_{\Rn} A(\xi) \ddd\xi,
 \end{equation}
where $$A(\xi) f(x):= \int_{\Rn} c(x,y,\xi)\, e^{i\varphi(x,\xi)-i\varphi(y,\xi)}\, f(y) \dd y.$$
Noting that $A(\xi)=0$ for $\xi$ outside some compact set, we observe that condition $(i)$ of Lemma \ref{calderon-vaillancourt lemma} follows from Young's inequality and \eqref{derivative estimates for c} with $\sigma=0,$  and condition $(ii)$ of Lemma \ref{calderon-vaillancourt lemma} follows from the assumption of the compact support of the amplitude. To verify condition $(iii)$ we confine ourselves to the estimate of $\Vert A^{\ast}(\xi_1) A(\xi_2)\Vert$, since the one for $\Vert A(\xi_1) A^{\ast}(\xi_2)\Vert$ is similar. To this end,
a calculation shows that the kernel of $A^{\ast}(\xi_1) A(\xi_2) $ is given by
\begin{equation}\label{kernel of A star A}
K(x,y,\xi_1 , \xi_2):= \int_{\Rn} \overline{c}(z,x,\xi_1)\, c(z,y,\xi_2 )\, e^{i[\varphi(z,\xi_{2})-\varphi(z,\xi_{1})+\varphi(x,\xi_{1})-\varphi(y,\xi_{2})]} \dd z.
\end{equation}
The estimate \eqref{derivative estimates for c} yields
\begin{equation}\label{first estimate for K}
|K(x,y,\xi_1 , \xi_2)|\lesssim \langle \xi_1\rangle ^{m_1}\, \langle \xi_2\rangle ^{m_1} \int_{\Rn}   \Lambda(\langle \xi_{1}\rangle ^{\rho} (x-z))\,
\Lambda(\langle \xi_{2}\rangle ^{\rho} (y-z)) \dd z.
\end{equation}
Therefore by choosing $N$ large enough, Young's inequality and using the fact that $\int_{\Rn} \Lambda (x) \, dx \lesssim 1$ yield
\begin{equation}\label{first estim for A star A}
\Vert A^{\ast}(\xi_1) A(\xi_2)\Vert \lesssim \langle \xi_1\rangle ^{m_1 -n\rho}\, \langle \xi_2\rangle ^{m_1 -n\rho}.
\end{equation}
At this point we introduce another first order differential operator $M:= G^{-2} \{1-i(\langle\nabla_{z}\varphi(z,\xi_2)-\nabla_{z}\varphi(z,\xi_1),\nabla_{z}\rangle)\}$, with $G=(1+|\nabla_{z}\varphi(z,\xi_2)-\nabla_{z}\varphi(z,\xi_1)|^2)^{\frac{1}{2}}.$ Using the fact that $M e^{i(\varphi(z,\xi_{2})-\varphi(z,\xi_{1}))}= e^{i(\varphi(z,\xi_{2})-\varphi(z,\xi_{1}))},$ integration by parts in \eqref{kernel of A star A} yields
\begin{equation}
    \int_{\Rn} ({}^{t} M)^{N'}\{\overline{c}(z,x,\xi_1)\, c(z,y,\xi_2 )\}\, e^{i[\varphi(z,\xi_{2})-\varphi(z,\xi_{1})+\varphi(x,\xi_{1})-\varphi(y,\xi_{2})]} \dd z.
\end{equation}
Using the second part of Lemma \ref{integration by parts lem}, we find that $({}^{t} M)^{N'}\{\overline{c}(z,x,\xi_1)\, c(z,y,\xi_2 )\}$ is a linear combination of terms of the form
\begin{equation}\label{differential operator M}
    G^{-k}\bigg\{\prod_{\nu =1}^{q}(\partial^{\beta_{\nu}}_{z}\varphi(z,\xi_2)-\partial^{\beta_{\nu}}_{z}\varphi(z,\xi_1))\bigg\} \partial^{\gamma_1}_{z}\overline{c}(z,x,\xi_1)\, \partial^{\gamma_2}_{z} c(z,y,\xi_2),
\end{equation}
where $k,$ $q,$ $\beta_\nu$ satisfy the inequalities in \eqref{relations for order of derivatives} and $|\gamma_1|+|\gamma_2| \leq N'.$\\

Now we observe that \eqref{lowerbound for phi in xi} yields that

$$G^{-k}\lesssim (1+|\xi_1 -\xi_2|)^{-k}.$$

 Moreover using \eqref{eq:L2 condition_new} we can also deduce that 
 $$|\partial^{\beta_{\nu}}_{z}\varphi(z,\xi_2)-\partial^{\beta_{\nu}}_{z}\varphi(z,\xi_1)|\lesssim |\xi_1-\xi_2|\leq 1+|\xi_1-\xi_2|.$$ Moreover  \eqref{relations for order of derivatives} yields that $q-k\leq -N'$ and therefore we obtain
 $$ G^{-k}\Big|\prod_{\nu =1}^{q}(\partial^{\beta_{\nu}}_{z}\varphi(z,\xi_2)-\partial^{\beta_{\nu}}_{z}\varphi(z,\xi_1))\Big|\lesssim (1+|\xi_1-\xi_2|)^{q-k}\lesssim  (1+|\xi_1-\xi_2|)^{-N'} \lesssim |\xi_1- \xi_2|^{-N'}. $$

Thus \eqref{derivative estimates for c} and \eqref{differential operator M}, yield the following estimate for $K(x,y,\xi_1 , \xi_2)$
\begin{align}\label{second estimate for K}
 &|K(x,y,\xi_1 , \xi_2)|\lesssim |\xi_1 -\xi_2|^{-N'} \langle \xi_1\rangle ^{m_1}\, \langle \xi_2\rangle ^{m_1}(1+|\xi_1|+|\xi_2|)^{\delta N'}  \\ \nonumber
 &\qquad\times \int_{\Rn}  \Lambda(\langle \xi_{1}\rangle ^{\rho} (x-z))\, \Lambda(\langle \xi_{2}\rangle ^{\rho} (y-z)) \dd z.
\end{align}
Once again, choosing $N$ large enough, Young's inequality yields
\begin{equation}\label{second estim for A star A}
\Vert A^{\ast}(\xi_1) A(\xi_2)\Vert \lesssim \langle \xi_1\rangle ^{m_1-n\rho}\, \langle \xi_2\rangle ^{m_1-n\rho} \frac{(1+|\xi_1|+|\xi_2|)^{\delta N'}}{|\xi_1 -\xi_2|^{N'}}.
\end{equation}
Using the fact that for $x>0,$ $\inf (1,x) \sim (1+\frac{1}{x})^{-1}$, one optimizes the estimates \eqref{first estim for A star A} and \eqref{second estim for A star A} by
\begin{align}\label{optimal estim for A star A}
\Vert A^{\ast}(\xi_1) A(\xi_2)\Vert &\lesssim \langle \xi_1\rangle ^{m_1 -n\rho}\, \langle \xi_2\rangle ^{m_1 -n\rho}
\bigg(1+ \frac{|\xi_1 -\xi_2|^{N'}}{(1+|\xi_1|+|\xi_2|)^{\delta N'}}\bigg)^{-1}
\\\nonumber &:= h^{2}(\xi_1 , \xi_2).
\end{align}
Therefore  recalling that $m_1=n(\rho-\delta),$ in applying Lemma \ref{calderon-vaillancourt lemma}, we need to show that
\begin{equation}
  K(\xi_1 , \xi_2)= (1+|\xi_1|)^{\frac{-n\delta}{2}} (1+|\xi_2|)^{-\frac{n\delta}{2}}  \bigg(1+ \frac{|\xi_1 -\xi_2|^{N'}}{(1+|\xi_1|+|\xi_2|)^{\delta N'}}\bigg)^{-\frac{1}{2}}
\end{equation}
is the kernel of a bounded operator in $L^2.$ At this point we use Schur's lemma, which yields the desired conclusion provided that
$$ \sup_{\xi_1} \int_{\Rn} K(\xi_1 , \xi_2) \dd \xi_2, \quad \sup_{\xi_2} \int_{\Rn} K(\xi_1 , \xi_2) \dd \xi_1 $$
are both finite. Due to the symmetry of the kernel, we only need to show the finiteness of one of these quantities.\\
\noindent To this end, we fix $\xi_1$ and consider the domains $\mathcal{A}=\{(\xi_1, \xi_2);\, |\xi_2| \geq 2 |\xi_1|\},$ $\mathcal{B}=\{(\xi_1, \xi_2);\, \frac{|\xi_1|}{2}\leq |\xi_2| \leq 2 |\xi_1|\},$ and $\mathcal{C}=\{(\xi_1, \xi_2);\, |\xi_2| \leq \frac{ |\xi_1|}{2}\}.$ Now we observe that on the set $\mathcal{A},$ $K(\xi_1, \xi_2)$ is dominated by
\begin{equation*}
  (1+|\xi_1|)^{-\frac{n\delta}{2}} (1+|\xi_2|)^{-\frac{n\delta}{2}+\frac{N'}{2}(\delta -1)},
\end{equation*}
on $\mathcal{B},$ $K(\xi_1, \xi_2)$ is dominated by
\begin{equation*}
  (1+|\xi_1|)^{-n \delta}  \bigg(1+ \frac{|\xi_1 -\xi_2|^{N'}}{(1+|\xi_1|)^{\delta N'}}\bigg)^{-\frac{1}{2}},
\end{equation*}
and on $\mathcal{C},$ $K(\xi_1, \xi_2)$ is dominated by
\begin{equation*}
   (1+|\xi_2|)^{-\frac{n\delta}{2}} (1+|\xi_1|)^{-\frac{n\delta}{2}+\frac{N'}{2}(\delta -1)}.
\end{equation*}
Therefore, if $I_{\Omega}:= \int_{\Omega} K(\xi_1 , \xi_2 ) \ddd \xi_2,$ then choosing $\frac{N'}{2}(\delta -1)<-n ,$ which is only possible if $\delta<1$, we have that $I_{\mathcal{A}} <\infty$ uniformly in $\xi_1$. Also,
\begin{equation}
  I_{\mathcal{C}} \leq (1+|\xi_1|)^{n-\frac{n\delta}{2}+\frac{N'}{2}(\delta -1)}\leq C,
\end{equation}
which is again possible by the fact that $\delta<1$ and a suitable choice of $N'.$
In $I_\mathcal{B}$ let us make a change of variables to set $\xi_2 -\xi_1 = (1+|\xi_1| )^{\delta} \eta$, then
\begin{equation}
I_{\mathcal{B}} \leq \int_{\Rn} (1+|\eta|^{N'}) ^{-\frac{1}{2}} \ddd \eta <\infty,
\end{equation}
by taking $N'$ large enough. These estimates yield the desired result and the proof of their theorem is therefore complete.
\end{proof}

\section{A \texorpdfstring{$h^p\to L^p$}{} estimate for oscillatory integral operators}\label{sec:hplp}

In this section, we show regularity results for oscillatory integral operators. Apart from the distinction caused by the type of the amplitudes, the values of $k$ also play a decisive role in their regularity theory. When $0<k<1$, it turns out that the type of the amplitude i.e. $\rho$ and $\delta$ can be incorporated in the analysis in such a way that the critical order of decay of the amplitude can be improved compared to the case of $k\geq 1$ where the order of the amplitude is
\begin{equation}\label{decay1}
    -kn\Big|\frac{1}{p}-\frac{1}{2}\Big|+\min\Big\{0,\frac{n}{2}(\rho-\delta)\Big\}.
\end{equation}
To this end, let \begin{equation}\label{def kappa}\varkappa=\min(\rho,1-k)\end{equation} 
and set 
\begin{equation}\label{decay}
    m(p):= m_\varkappa(p)+ \zeta.
\end{equation}
where
$$m_\varkappa(p):=-n(1-\varkappa)\Big|\frac{1}{p}-\frac{1}{2}\Big|$$

and
$$\zeta:=\min\Big\{0,\frac{n}{2}(\rho-\delta)\Big\}.$$
Observe that for $k\geq 1$ one recovers \eqref{decay1} from \eqref{decay} since then $\varkappa=1-k$.\\

\begin{Lem}\label{L2 of kernel lemma}
Let $0\leq \rho\leq 1$ and $0\leq \delta< 1$. Assume that $\varphi\in \textart F^k$ is \emph{SND} and satisfies the $L^2$-condition \eqref{eq:L2 condition_old} for $0<k<\infty$. 
Then if $a(x,\xi)\in S^{m}_{\rho,\delta}(\Rl^n)$ for some $m\in\Rn$, let

\begin{equation*}
    K_{j}(x,y) = \int_{\Rl^n} e^{i(\varphi(x,\xi)-y\cdot \xi)}\,\sigma_{j}(x,\xi)  \ddd \xi.
\end{equation*}
where $\sigma_{j}(x,\xi)=\psi_j(\xi)\,a(x,\xi)$. For $M\geq 0$, $y\in \Rn$ and $j\geq 1$ we have
\begin{align}\label{L2 uppskattning karnan}
    \Vert (1+2^{j\varkappa} |x-y|)^M \partial^{\beta}_y K_{j}(x,y)\Vert_{L^2_x(\Rl^n)}
&\lesssim 2^{j(\tilde m+n/2+\vert \beta\vert)},
\end{align}
where $\tilde m:=m-\zeta$. 
Moreover, under the extra assumption of $\varphi(x,0)=0$
then estimate \eqref{L2 uppskattning karnan} is also valid when $K_j(x,y)$ is replaced by the kernel of the adjoint $K_{j}^{*}(x,y)$.
\end{Lem}

\begin{proof}
Observe that for any multi-index $\beta$ and any $j\geq 0$ we have
\begin{align*}
\partial^{\beta}_{y}K_{j}(x,y) 
& = \int_{\Rl^n}  \psi_j(\xi)\,
\partial^{\beta}_{y} e^{i(\varphi(x,\xi )-y\cdot \xi)} \,a(x,\xi) \ddd \xi \\
& = 
\int_{\Rl^n}
e^{i(\varphi(x,\xi )-y\cdot \xi)}\, \sigma_j^{\beta}(x,\xi) \ddd\xi,
\end{align*}

where
$$\sigma_{j}^{\beta}(x,\xi)
:=  \sigma_{j}(x,\xi)\,(-i\xi)^\beta.
$$
Therefore, since $|\partial^{\alpha}_\xi \partial^{\gamma}_x\sigma_{j}^{\beta}(x,\xi)|\lesssim 2^{j(m+|\beta|-\rho|\alpha|+\delta|\gamma|)}$ we can reduce ourselves to the case when $\beta=0$.\\

Now, using \eqref{amplitude derivative estim} and \eqref{linearizeddecayestimate} and letting  $\Psi_j$ be a Littlewood-Paley function that is supported on a larger annulus in the sense of Definition \ref{def:LP}, we have
\begin{align*}
&(x-y)^\alpha K_j(x,y)
=\int_{\Rl^n} (-i)^{\vert \alpha\vert}\partial_\xi^\alpha\,e^{i(x-y)\cdot \xi}\,\sigma_{j}(x,\xi)\, e^{i(\varphi(x,\xi)-x\cdot \xi)}\ddd\xi\\
& \quad =\int_{\Rl^n} i^{\vert \alpha\vert}\partial_\xi^\alpha \Big[\,e^{i(\varphi(x,\xi)-x\cdot \xi)}\,\sigma_{j}(x,\xi)\Big]\, e^{i(x-y)\cdot \xi}\,{\Psi}_j (\xi)\ddd\xi \\
& \quad = \sum_{\alpha_1 + \alpha_2=\alpha} \!\!C_{\alpha_1, \alpha_2}\int_{\Rl^n}
\partial_\xi^{\alpha_1} \sigma_j(x,\xi) \,
\partial_\xi^{\alpha_2}e^{i(\varphi(x,\xi)-x\cdot \xi)}\,\,e^{i(x-y)\cdot \xi} {\Psi}_j (\xi)\ddd\xi \\
& \quad =\!\! \sum_{\substack{\alpha_1 + \alpha_2=\alpha\\ \lambda_1 + \dots + \lambda_r = \alpha_2} }\!\!\!\!
C_{\alpha_1, \alpha_2,\lambda_1 ,\dots , \lambda_r} 
2^{j(\tilde m-\varkappa|\alpha|)}   \!\!
\int_{\Rl^n} b_j^{\alpha_1,\alpha_2, \lambda_1, \dots \lambda_r}(x,\xi)\,
e^{i\varphi(x,\xi)}\, e^{-iy\cdot\xi}\,{\Psi}_j (\xi) \ddd\xi \\
& \quad =:\!\! \sum_{\substack{\alpha_1 + \alpha_2=\alpha \\ \lambda_1 + \dots + \lambda_r = \alpha_2} }
\!\!\!\!C_{\alpha_1, \alpha_2, \lambda_1, \dots, \lambda_r} 
2^{j(\tilde m-\varkappa|\alpha|)}  
S_{j}^{\alpha_1, \alpha_2, \lambda_1, \dots \lambda_r}\big(\tau_{-y}\widehat\Psi_j\big)(x),
\end{align*}

where $\tau_{-y}$ is a translation by $-y$, $|\lambda_j| \geq 1$ and
\begin{align*}
b_j^{\alpha_1,\alpha_2, \lambda_1, \dots, \lambda_r}(x,\xi)
& := 2^{-j(\tilde m-\varkappa|\alpha|)}\\
& \qquad \times \partial_\xi^{\alpha_1} \sigma_j(x,\xi) \, 
\partial_\xi^{\lambda_1}(\varphi(x,\xi)-x\cdot \xi)
\dots
\partial_\xi^{\lambda_r}(\varphi(x,\xi)-x\cdot \xi),
\end{align*}
\\

Now we claim that 
$b_j^{\alpha_1,\alpha_2, \lambda_1, \dots, \lambda_r}\in S^{\zeta}_{0,\delta}(\Rl^n),$ uniformly in $j$. Indeed, since \linebreak$a \in S^{m}_{\rho,\delta}(\Rl^n)$ and $\varphi\in \textart F^k$, we can write
\begin{align*}
|b_j^{\alpha_1,\alpha_2,\alpha_3, \lambda_1, \dots, \lambda_r}(x,\xi)|
&\lesssim 2^{-j(\tilde m-\varkappa|\alpha|)} \, 
2^{j\tilde m+j\zeta-j\rho|\alpha_1|} \,
2^{j(k-1)|\alpha_2|} \\
&\lesssim 
2^{j\zeta}2^{j\varkappa|\alpha|} \, 
2^{-j\rho|\alpha_1|} \,
2^{j(k-1)|\alpha_2|}\\
&\lesssim 2^{j\zeta}.
\end{align*}

Moreover, observe also that by the $\textart F^k$-condition for $k<1$,
\begin{equation*}
|\partial_x^{\beta_{j}}\partial_\xi^{\lambda_{j-1}+\gamma_{j}}(\varphi(x,\xi)-x\cdot \xi)|
    \lesssim 2^{j(k-|\lambda_{j-1}+\gamma_{j}|)}  \text{ when}\,\,\,|\beta_{j}| \geq 0
\end{equation*}
and for $k\geq 1$ the $L^2$-condition \eqref{eq:L2 condition_old}  yields that 
\begin{equation*}
|\partial_x^{\beta_{j}}\partial_\xi^{\lambda_{j-1}+\gamma_{j}}(\varphi(x,\xi)-x\cdot \xi)|
    \lesssim 1\text{ when}\,\,\,|\beta_{j}| \geq 1.
\end{equation*}
Thus, using these estimates, we have for all $k<0$ and all $1\leq j\leq r+1$ that
\begin{equation}
\label{eq:mixed_derivatives_linear_phase}
|\partial_x^{\beta_{j}}\partial_\xi^{\lambda_{j-1}+\gamma_{j}}(\varphi(x,\xi)-x\cdot \xi)|
    \lesssim 2^{j(k-1)|\lambda_{j-1}+\gamma_{j}|},  \text{ when}\,\,\,|\beta_{j}| \geq 1
\end{equation}
Now using \eqref{eq:mixed_derivatives_linear_phase} we can also check that, for any multi-indices $\gamma$ and $\beta$, 
\begin{align}\label{eq:bjs}
    &|\partial_{\xi}^\gamma \partial_x^\beta \, b_j^{\alpha_1,\alpha_2, \lambda_1, \dots, \lambda_r}(x,\xi)|\\
    &\nonumber\lesssim\sum_{\substack{\gamma_1 + \dots + \gamma_{r+1} = \gamma \\ \beta_1 + \dots + \beta_{r+1} = \beta} }2^{-j(\tilde m-\varkappa|\alpha|)} |\partial_x^{\beta_{1}}\partial_\xi^{\alpha_1+\gamma_1} \sigma_j(x,\xi)| \, 
    |\partial_x^{\beta_{2}}\partial_\xi^{\lambda_1+\gamma_2}(\varphi(x,\xi)-x\cdot \xi)|\dots
    \\
    &\nonumber \qquad \times
    |\partial_x^{\beta_{r+1}}\partial_\xi^{\lambda_r+\gamma_{r+1}}(\varphi(x,\xi)-x\cdot \xi)|\\
    &\nonumber\lesssim\sum_{\substack{\gamma_1 + \dots + \gamma_{r+1} = \gamma \\ \beta_1 + \dots + \beta_{r+1} = \beta} }2^{-j(\tilde m-\varkappa|\alpha|)}
     2^{j(\tilde m+\zeta+ \delta|\beta_1|-\rho|\alpha_1+\gamma_1|)}2^{j(k -1)|\alpha_2|}\\
    &\nonumber\lesssim\sum_{\substack{\gamma_1 + \dots + \gamma_{r+1} = \gamma \\ \beta_1 + \dots + \beta_{r+1} = \beta} }2^{-j(\tilde m-\varkappa|\alpha|)}
     2^{j(\tilde m+\zeta+\delta|\beta_1|)} 2^{-j\varkappa|\alpha|}\\
    &\nonumber\lesssim 2^{j(\zeta+\delta|\beta|)},
\end{align}

hence $b_j^{\alpha_1,\alpha_2, \lambda_1, \dots, \lambda_r}\in S^\zeta_{0,\delta}(\Rl^n).$  \\

Therefore, Theorem \ref{thm:Calderon-Vaillancourt for OIOs} yields that 
\begin{align*}
&\Vert (x-y)^\alpha K_{j}(x,y)\Vert_{L^2_x(\Rl^n)}\\
& \lesssim \!\!\sum_{\substack{\alpha_1 + \alpha_2=\alpha \\ \lambda_1 + \dots + \lambda_r = \alpha_2} }\!\! 2^{j(\tilde m-\varkappa|\alpha|)} 
\big\|S_{j}^{\alpha_1, \alpha_2, \lambda_1, \dots \lambda_r}\big(\tau_{y}\widehat\Psi_j\big)\big\|_{L^2(\Rl^n)} \\
& \lesssim 2^{j(-\vert\alpha\vert\varkappa +\tilde m)}\big\|\widehat\Psi_j\big\|_{L^2(\Rl^n)} \\
& \lesssim 2^{j(-\vert\alpha\vert \varkappa+\tilde m+n/2)}.
\end{align*}
From this and the discussion at the beginning of the proof one can deduce that
\begin{equation}
\Vert (x-y)^\alpha \partial_{y}^{\beta}K_{j}(x,y)\Vert_{L^2_x(\Rl^n)}\\
\lesssim 2^{j(-\vert\alpha\vert\varkappa+|\beta|+\tilde m+n/2)}.
\end{equation}
Thus by summing over all $|\alpha|\leq M$ for any integer $M$ one obtains 
\begin{align*}
\Vert (1+2^{j\varkappa} |x-y|)^M \partial_{y}^{\beta}K_{j}(x,y)\Vert_{L^2_x(\Rl^n)}
\lesssim 2^{j(\tilde m+n/2+|\beta|)}.
\end{align*}
For the kernel of the adjoint, we have that 
\begin{equation*}
K_{j}^{*}(x,y) = \int_{\Rl^n}  \,e^{-i(\varphi(y,\xi )-x\cdot \xi)}\,\overline{\sigma_j (y,\xi)}  \ddd \xi,
\end{equation*}
therefore for any multi-index $\alpha$ we have
\begin{align*}
(x-y)^\alpha\partial^{\beta}_{y}K_{j}^{*}(x,y) 
& = (x-y)^\alpha\int_{\Rl^n}  \,
\partial^{\beta}_{y} \Big(e^{-i(\varphi(y,\xi )-x\cdot \xi)}\,\overline{\sigma_j (y,\xi)}\Big) \ddd \xi \\
& =(x-y)^\alpha 
\int_{\Rl^n} e^{-i(\varphi(y,\xi )-x\cdot \xi)}\,
\sigma_{j}^{\beta,*}(y,\xi) \ddd\xi\\
& = 
\int_{\Rl^n} e^{i(x-y)\cdot\xi}\,
i^{\vert \alpha\vert}\partial_\xi^\alpha \Big[e^{-i(\varphi(y,\xi )-y\cdot \xi)}\,
\sigma_{j}^{\beta,*}(y,\xi)\Big] \ddd\xi,
\end{align*}
where
$$\sigma_{j}^{\beta,*}(y,\xi)
:=  \sum_{\substack{\beta_1 + \beta_2=\beta \\ \lambda_1 + \dots + \lambda_r = \beta_2} }
C_{\beta_1, \beta_2, \lambda_1, \dots \lambda_r}\,
\partial_y^{\beta_1} \overline{\sigma_j(y,\xi)}\,
\partial_y^{\lambda_1}\varphi(y,\xi)
\cdots
\partial_y^{\lambda_r}\varphi(y,\xi),
$$
and $|\lambda_j| \geq 1$.
Now, for $|\lambda_j + \beta| \geq 2$,
$$|\partial^{\lambda_j}_y \partial^{\beta}_{\xi} \varphi(y,\xi)|
\lesssim 1,$$ 
and using that $\varphi(y,0)=0$, the mean-value theorem and $L^2$-condition \eqref{eq:L2 condition_old}, we obtain $$| \nabla_y \varphi(y,\xi)|
\lesssim |\xi|.$$ 
From these estimates, we deduce that for any multi-index $\gamma$ one has
$|\partial^{\alpha}_y\partial^{\gamma}_\xi\sigma_{j}^{\beta,*}(y,\xi)|\lesssim 2^{j(m+|\beta|-\rho|\gamma|+\delta\vert\alpha\vert))}$.
Therefore, following the same line of reasoning as in the case of $K_j(x,y)$ yields the estimate given in \eqref{L2 uppskattning karnan} for $K_{j}^{*}(x,y)$. 
\end{proof}

Now we are ready to show the main result of this section.

\begin{Th}\label{thm:hp-Lp_oio}
Let $n\geq 1$, $0<k<\infty$, and $0<p<\infty$. Assume that $\varphi\in \textart F^k$ {is \emph{SND}}, satisfies the \emph{LF}$(\mu)$-condition for some $0<\mu\leq 1$, and the $L^2$-condition \eqref{eq:L2 condition_old}. Assume also that $a(x,\xi)\in S^{m(p)}_{\rho,\delta}(\Rl^n),$ for $0\leq \rho\leq 1$ and $0\leq \delta< 1$. Then the \emph{OIO} $T_a^\varphi$ is bounded from
\begin{enumerate}
    \item[$(i)$] $h^{p}(\Rl^n)\to L^p(\Rl^n),$ and
    \item[$(ii)$] $L^\infty(\Rl^n)\to \bmo(\Rl^n)$ provided that one also has $|\nabla \varphi_x (x,0)|\in L^{\infty}(\Rl ^n)$. 
\end{enumerate}
\end{Th}

\begin{proof}
Let $\chi\in C_c^\infty(\Rn)$ be supported in $\{\xi:|\xi|\lesssim 1\}$, and write
    \begin{align*}
        T_a^\varphi f(x) &: =
        \int_{\Rn} e^{i\varphi(x,\xi)} a(x,\xi) \widehat{f}(\xi) (1-\chi(\xi)) \ddd\xi 
        +
        \int_{\Rn} e^{i\varphi(x,\xi)} a(x,\xi) \widehat{f}(\xi) \chi(\xi) \ddd\xi\\
        &=T_{\text{high}} f(x)+T_{\text{low}} f(x).
    \end{align*}
    The boundedness of $T_{\text{low}}$ follows from the low frequency result Theorem \ref{thm:low_freq_TL_BL_OIO}. Thus for the remainder of the proof we only consider the high frequency portion of the operator.\\

Let $T_j^\nu$ be the operator associated to the kernel in Definition \ref{second decomposition}, so that 
\begin{equation*}
    T_{\text{high}}=\sum_{j=1}^\infty \sum_{\nu=1}^{\mathscr{N}_j} T_j^\nu.
\end{equation*}

\textbf{Case when $\boldsymbol{0<p<1}$:\\}

First we consider the case when $0<p<1$. Let $\mathfrak{a}$ be a $p$-atom supported in a cube $Q$ with side length $l_Q$ and let $2Q$ be the cube with the same center and twice the side length. Since $0<p<1$ we have
\begin{align}\label{eq:twoterms2}
\Vert T_{\text{high}} \at\Vert^p_{L^p(\Rn)} 
& \lesssim \Vert T_{\text{high}}\at\Vert^p_{L^p(2Q)} +\sum_{j=1}^\infty 
\Vert T_j \at\Vert^p_{L^p(\Rn \setminus 2Q)} 
\end{align}

Observe that by H\"older's inequality and the $L^2$-boundedness we have,
\begin{align}\label{eq:termauna}
\Vert T_{\text{high}}\at\Vert^p_{L^p(2Q)} 
&\lesssim \Vert T_{\text{high}}\at\Vert^2_{L^2(2Q)} \Vert 1\Vert^p_{L^{2p/(2-p)}(2Q)}
\lesssim \Vert\at\Vert^2_{L^2(2Q)} l_Q^{n(2-p)/2p}\\
&\nonumber\lesssim l_Q^{-n(2-p)/2p} l_Q^{n(2-p)/2p}\lesssim 1
\end{align}

By Lemma \ref{L2 of kernel lemma} we have for $0<k<1$ that 
\begin{align*}
\Vert (1+2^{-j\varkappa} |x-y|)^M K_{j}(x,y)\Vert_{L^2_x(\Rl^n)}
\lesssim 2^{j(m_{\varkappa}(p)+n/2)},
\end{align*}
Observe that for $t\in[0,1]$ and $x\in \Rl^n\setminus 2Q$, one has 
\begin{equation}\label{usefulinequality}
    |x-\bar{y}|\lesssim |x-y+t(y-\bar{y})|
\end{equation}
Now, setting  
$$g(x)=\jap{2^{j\varkappa}|x-\bar y|}^{-M}$$
Observe that we have for $q\geq 1$,
\begin{equation}\label{metric norm1}
    \|g\|_{L^\tau(\Rl^n)}\lesssim 2^{nj\varkappa/\tau}.
\end{equation}

By H\"older's inequality and Lemma \ref{L2 of kernel lemma} and \eqref{metric norm1} we have for $l_Q>1$ and $w\geq 1$ that
\begin{align}\label{comple}
    \Vert T_j \at\Vert_{L^p(\Rn\setminus 2Q)}
    & \lesssim \Big\| \frac{1}{g(x)} \int_Q  K_j (x,y)\at(y) \dd y\Big\|_{L^2(\Rl^n)}\|g\|_{L^\tau(\Rl^n)}\\
    &\nonumber \lesssim  \int_Q\Big\| \frac{1}{g(x)} K_j (x,y)\Big\|_{L^2(\Rl^n)} |\at(y)| \dd y\,\|g\|_{L^\tau(\Rl^n)}\\
    &\nonumber \lesssim  2^{-nj\varkappa/\tau} \int_Q\Big\| \frac{1}{g(x)}K_j (x,y)\Big\|_{L^2(\Rl^n)} |\at(y)| \dd y\\
    &\nonumber \lesssim  l_Q^{n-n/p}2^{j(m_{\varkappa}(p)+n/2-n\varkappa(\frac{1}{p}-\frac{1}{2}))} \lesssim  l_Q^{n-n/p} 2^{j(n-n/p)}
\end{align}
where $\frac{1}{\tau}=\frac{1}{p}-\frac{1}{2}$.\\

Now, if $l_Q<1$, Taylor expansion of $K$ in the $y$-variable around $\bar y$, using the moment conditions of $\at$ yields for all $t\in[0,1]$ and $N:=[n(1/p-1)]$ yield that
\begin{align*}
    &\|T_j \at\|_{L^p(\Rn \setminus 2Q)}^p 
    \lesssim \sum_{|\alpha|=N+1}
    \int_{\Rn \setminus 2Q} 
    \Big(\int_Q |{ \partial_y^\alpha K_j(x,t(y-\bar y))}||y-\bar y|^{N+1}|\at(y)|\dd y \Big)^p \dd x.\nonumber
\end{align*}
Moreover, by using that $|y-\bar y|^{N+1}\lesssim r^{N+1}$ with \eqref{usefulinequality} and Lemma \ref{L2 of kernel lemma} we obtain the following estimate by a similar calculation as in \eqref{comple},
\begin{align}\label{complementnorm of taylor expansion1}
    &\|T_j \at\|_{L^p(\Rn \setminus 2Q)}^p 
     \lesssim 2^{jm_{\varkappa}(p)p+jnp/2-jnp\varkappa(1/p-1/2)}2^{jp(N+1)}l_Q^{np-n+p(N+1)}.
\end{align}
Now, since $l_Q<1$, take the unique integer $k_Q\in \mathbb Z_{>1}$ such that $2^{-k_Q}< l_Q\leq 2^{-(k_Q-1)}$. Then using \eqref{comple} and \eqref{complementnorm of taylor expansion1} we have
\begin{align*}
    &\sum_{j=1}^\infty\Vert T_j \at\Vert^p_{L^p(\Rn \setminus 2Q)} \\ 
    &\nonumber\qquad \lesssim \sum_{j= k_Q}^\infty 2^{jp(m_{\varkappa}(p)+n/2+n\varkappa(\frac{1}{p}-\frac{1}{2}))}l_Q^{np-n}\\ 
    &\nonumber\qquad\qquad+ \sum_{j=1}^{k_Q-1} 2^{jm_{\varkappa}(p)p+np/2+jnp\varkappa(1/p-1/2)}2^{jp(N+1)}l_Q^{np-n+p(N+1)}\\ 
    &\nonumber\qquad \lesssim \sum_{j= k_Q}^\infty 2^{j(np-n)}l_Q^{np-n}+ \sum_{j=1}^{k_Q-1} 2^{j(np-n+p(N+1))}l_Q^{np-n+p(N+1)}\\ 
    &\nonumber\qquad \lesssim  2^{k_Q(np-n)}l_Q^{np-n} + 2^{k_Q(np-n+p(N+1))}l_Q^{np-n+p(N+1)}\\ 
    &\nonumber\qquad \lesssim  l_Q^{-(np-n)}l_Q^{np-n} + l_Q^{-(np-n+p(N+1))}l_Q^{np-n+p(N+1)}\\ 
    &\nonumber\qquad\sim 1.
\end{align*}

Now, if $l_Q\geq 1$ we do the same calculation as above but with $k_Q=0$, and do not consider the case $j<k_Q$. Thus we conclude that
\begin{align}\label{eq:twoterms2proved}
\Vert T_{\text{high}} \at\Vert^p_{L^p(\Rn)} 
& \lesssim \Vert T_{\text{high}}\at\Vert^p_{L^p(2Q)} +\sum_{j=1}^\infty 
\Vert T_j \at\Vert^p_{L^p(\Rn \setminus 2Q)} \lesssim 1.
\end{align}

\textbf{Interpolation and arguments for the adjoint operator:\\}

Now, interpolating the result abovewith the $L^2$-boundedness result in Theorem \ref{thm:Calderon-Vaillancourt for OIOs} for $0\leq \rho\leq 1$ and $0\leq\delta<1$ using Riesz-Thorin's interpolation theorem, one obtains the $h^p-L^p$--boundedness of $T_a^\varphi$ (with the decay $m(p)$), in the range $0<p\leq 2$.\\

Next we observe that one may write the phase as $\varphi(x, \xi)=\varphi(x,\xi)-\varphi(x, 0)+ \varphi(x,0)=:\psi(x,\xi)+\varphi(x,0)$. This would certainly yield that $$\|T^\varphi_a f\|_{L^p(\Rl^n)} = \|T^\psi_a f\|_{L^p(\Rl^n)},$$ for $2\leq p<\infty$. Moreover observe that $\psi(x,0)=0$, and therefore we can without loss of generality assume that $\varphi(x,0)=0$ in $T_a^\varphi$. Now in order to prove the $h^p-L^p$--boundedness of $T_a^\varphi$ for $2\leq p<\infty$, using duality and interpolation, it is enough to show that the adjoint operator ${(T_a^{\varphi})}^*$ is bounded from $h^p(\Rl^n)$ to $L^p(\Rl^n)$, for $0<p\leq 2.$  However, this can be shown by the same argument as in the proof of the $h^p-L^p$--boundedness of ${(T_a^{\varphi})}$, replacing the $L^2$-inequality for the kernel with the corresponding inequality for the kernel of the adjoint (given in Lemma \ref{L2 of kernel lemma}).\\

Now for the boundedness of $T_a^\varphi$ from $L^\infty(\Rl^n)$ to $\bmo(\Rl^n)$ one can write \linebreak$T_a^\varphi= e^{i\varphi(x,0)} T_a^{{\psi}}$ with  ${\psi}(x,0)=0.$  Then given the assumption that $\varphi\in\textart{F}^k$, that $\varphi$ satisfies the $L^2$-condition of Definition \ref{L2 conditions}, and the extra assumption \linebreak$|\nabla_x\varphi(x,0)|\in L^{\infty}(\Rl^n)$ on the phase function, one can use \eqref{pointwise multiplier} to reduce matters to the boundedness of $T_a^{{\psi}}$. But the boundedness of $T_a^{{\psi}}$ from $L^\infty(\Rl^n)$ to $\bmo(\Rl^n)$ is a consequence of the boundedness of $(T_a^{{\psi}})^{\ast}$ from $h^1(\Rl^n)$ to $L^1(\Rl^n)$ which was achieved above, and therefore the proof of the theorem is concluded.

 \end{proof}
\section{Transference to Triebel-Lizorkin spaces for \texorpdfstring{$0<p<1$}{} and \texorpdfstring{$q\geq p$}{}} \label{subsec:OIO_lowP_transference}

This section is devoted to one-half of the process of going from $h^p\to L^p$ to Triebel-Lizorkin boundedness. In particular, we state and prove the global boundedness of oscillatory integral operators with classical and exotic amplitudes on Triebel-Lizorkin spaces $F^s_{p,q}$ for $0<p<2$ and $q\geq p$. To this end, we define a molecular representation of the Triebel-Lizorkin spaces. Similar to the atomic representation of the Hardy spaces, this can be used to prove boundedness results. \\

Before we show the main results of this section we recall a number of useful lemmas from \cite{ATW}, and as the proofs are very similar or exactly the same we leave out the proofs. In particular, this is with regards to Lemma \ref{lem:TLatomseq1} through Lemma \ref{lem:local and global nonendpoint TL}.

\begin{Def}[Notation]
Let $\mathcal D$ be the set of all dyadic cubes in $\Rn$ and define the following sets:
\begin{enumerate}
    \item[$(i)$] $\mathcal D_j :=\set{Q\in \mathcal D: k_Q= j} $
    \item[$(ii)$] $\mathcal D_+:= \set{Q\in \mathcal D: k_Q\geq0}$
    \item[$(iii)$] $\mathcal D_j(Q):= \{\tilde  Q\in \mathcal D_{j}:\tilde Q\subseteq Q\} $.
\end{enumerate}
\end{Def}

Observe that for $j<k_Q$, $D_j(Q)=\emptyset.$\\

We start with defining a space of sequences that is easier to handle than the Triebel-Lizorkin spaces themselves.

\begin{Def}\label{def:TLatomseq1}For a sequence of complex numbers $b = \set{b_Q}_{\substack{ Q\in\mathcal D\\
l(Q)\leq 1}}$
we define
\eq{
g^{s,q}(b)(x) := \brkt{\sum_{ \tilde Q\in\mathcal D_+}
\big(2^{k_{\tilde Q}(s+n/2)}|b_{\tilde Q}|\chi_{\tilde Q}(x)\big)^q}^{1/q}.
}
We say that $b\in f^s_{p,q}$ is
\eq{
\| b\|_{f^s_{p,q}} := \|g^{s,q}(b)\|_{L^p(\Rn)}<\infty.
}
\end{Def}
In what follows, take $\check\Psi^{\tilde Q}(x) := 2^{-nk_{\tilde Q}/2} \check\Psi_{k_{\tilde Q}}(x- c_{\tilde Q})$, where $\Psi_{k_{\tilde Q}}$ as given in Definition \ref{def:LP} is the usual Littlewood-Paley piece. The following Lemma is a corollary of \cite[Theorem II B]{FJ:phi-transform}.
\begin{Lem}\label{lem:TLatomseq1}
Suppose $0 < p < \infty,$ $ 0 < q \leq \infty$, $s \in \Rl.$ For any sequence $b= \set{b_{\tilde Q}}_{\tilde Q\in\mathcal D}$ of complex numbers satisfying
$\|b\|_{
f^s_{p,q}} < \infty,$ one has
\eq{
    f(x) := \sum_{\tilde Q\in \mathcal D_+}b_{\tilde Q}\check\Psi^{\tilde Q}(x)
}
belongs to $F^s_{p,q}(\Rn)$ and
\eq{
\|f\|_{F^s_{p,q}(\Rn)}\lesssim 
\|b\|_{f^s_{p,q}}.
}
\end{Lem}
\begin{proof}
 This follows by combining \cite[Theorem II A i]{FJ:phi-transform} and \cite[Theorem 1]{HanPalWeiss}.  It is worth noting that this result can also be applied to inhomogeneous Triebel-Lizorkin spaces, as discussed in Chapter 12 of \cite{FJ:discrete-transform}.
\end{proof}

We now discuss the converse of Lemma \ref{lem:TLatomseq1}, namely when a Triebel-Lizorkin function can be expressed in terms of molecules and so-called "$\infty$-atoms", here denoted $b_{\iota,\tilde  Q}$.
\begin{Lem}[\cite{ATW}]\label{lem:TLatomseq}
Suppose $0 < p \leq 1,$ $p \leq q \leq \infty$. Every $f\in F^0_{p,q}(\Rn)$ has an atomic decomposition 
\eq{
    f(x)=\sum_{\iota=1}^\infty \lambda_\iota  \sum_{ \tilde Q\in \mathcal D_+}b_{\iota,\tilde  Q}\check\Psi^{\tilde  Q}(x), \qquad b_{\iota,\tilde  Q}\in 
    f^0_{\infty,q},
}
where $ \set{b_{\iota,\tilde  Q}}_{\tilde Q\in\mathcal D_+}=:b_\iota$ satisfies 
\eq{
\|b_\iota\|_{f^0_{\infty,q}} \leq 2^{k_{\tilde Q}n/p}. 
}
Moreover
\eq{
\|f\|_{F^0_{p,q}(\Rn)}\sim\inf_{\set{\lambda_\iota}}\brkt{\sum_{\iota=1}^\infty |\lambda_\iota|^p}^{1/p}
}

\end{Lem}

Next, we define the analog of the Hardy space atoms, which will be used to prove the boundedness results of our OIO's.

\begin{Lem}[\cite{ATW}]\label{lem:TLatomseq4}
Let $Q \in \mathcal D$, $0<p\leq 1$, $0<q\leq\infty$, $j\geq 1$ and $b = \set{b_{\tilde Q}}_{\tilde Q\in\mathcal D_+}\in f^0_{\infty,q}$ with $\|b\|_{f^0_{\infty,q}} \leq 2^{k_Qn/p}$. Define
\nm{eq:RQatom}{
    R_{Q,j}(x):=\sum_{\tilde Q\in\mathcal D_j( Q)}b_{\tilde Q}\check\Psi^{\tilde Q}(x).
}
Then $R_{Q,j}(x)=0$ for $j< k_{\tilde Q}$. Moreover, for $0<p<\infty$,
\eq{
\|R_{Q,j}\|_{L^{2}(\Rn)}\lesssim 2^{k_{ Q}n(1/p-1/2)}.
}
\end{Lem}

Now we start with the lifting results for OIO's to Triebel-Lizorkin spaces. In order to do that, we need to introduce a partition of unity to $R_{Q,j}$. We estimate the pieces separately. The first estimate is given in Lemma \ref{lem:TLatomseq5}.

\begin{Lem}[\cite{ATW}]\label{lem:TLatomseq5}
Suppose that $0<p<1$ and $q\geq p$. Let $a (x,\xi)\in S^{m(p)}_{\rho, \delta}(\Rn)$ be supported outside a neighborhood of the origin and assume that $\varphi\in \textart F^k$ {is \emph{SND}}, satisfies the \emph{LF}$(\mu)$-condition for some $0<\mu\leq 1$, and satisfies the $L^2$-condition \eqref{eq:L2 condition_old}. Then
\nm{eq:elände3}{
\| T_a^\varphi \psi_j(D)  \chi_{\Rn \setminus2\sqrt n  Q}R_{Q,j}\|_{L^p(\Rn )} \lesssim 2^{n(k_{Q}-j)}.
}
for $j\geq k_{ Q}$.
\end{Lem}

\begin{Lem}[\cite{ATW}]\label{lem:local and global nonendpoint TL}
Let $a(x,\xi)\in S_{\rho,\delta}^{m}(\Rn)$ and that $T_a^\varphi$ is an oscillatory integral operator that is bounded from $F_{p,p}^s(\Rn)$ to $F_{p,p}^s(\Rn)$. Assume also that the phase function $\varphi$ satisfies the conditions of Theorem \ref{thm:left composition with pseudo}. Then $T_\sigma^\varphi$ is bounded from $F_{p,q}^{s}(\Rn)$ to $F_{p,q}^s(\Rn)$ for $\sigma(x,\xi)\in S_{\rho,\delta}^{m-\varepsilon}(\Rn)$ where $\varepsilon>0$ is arbitrary.
\end{Lem}
\begin{proof}
    This lemma follows similarly to \cite[Lemma 5.8]{ATW}, with some minor modifications. One replaces the $h^p\to L^p$ boundedness result in that paper (Proposition 5.1) with Theorem \ref{thm:hp-Lp_oio} in this paper, beyond this abstract modification the proof remains the same, and therefore one only substitutes the hypothesis on $T_a^\varphi$ in \cite[Lemma 5.8]{ATW}, by the hypothesis necessary for the $h^p\to L^p$ boundedness of $T_a^\varphi$ in Theorem \ref{thm:hp-Lp_oio}.
\end{proof}

Now we are prepared to show the main lifting results of this section.

\begin{Lem}\label{lem:bstarcestimate}
Let $T_a^\varphi$ be an $\mathrm{OIO}$ with an amplitude $a\in S^{m(p)}_{\rho,\delta}(\Rn)$ for $\rho\in[0,1],\,\delta\in[0,1)$. Assume $\varphi\in \textart F^k$ {is \emph{SND}} and satisfies the $L^2$-condition \eqref{eq:L2 condition_old}. Moreover $T_j:= T_a^\varphi\psi_j(D)$, where $\psi_j(D)$ is a Littlewood-Paley piece as in \emph{Definition \ref{def:LP}}.  Furthermore, suppose that $f$ is supported in a cube $Q$ with
$
\|f\|_{L^{ 1}(\Rn)}\lesssim 2^{k_Qn(1/p-1)}.
$
\begin{enumerate}
    \item[$(i)$] If $0<p<1$, then
    \begin{align}\label{eq:LPoperator2}
        &&\sum_{j=\max\{k_Q+1,1\}}^\infty\| T_j f\|_{L^p(\Rn \setminus 2Q)}^p\lesssim 1.
    \end{align}
    \item[$(ii)$] If $0<p\leq 1$ and $f$ is an $h^p$-atom, see \emph{Definition \ref{def:hpatom}}, then
    \begin{align}\label{eq:LPoperator1}
        &&\sum_{j=1}^{\max\{k_Q,0\}}\| T_j  f\|_{L^p(\Rn \setminus 2Q)}^p\lesssim 1.
    \end{align}
    
\end{enumerate}
Moreover, the same estimates hold true for the adjoint operator $(T_j)^*.$

\end{Lem}
\begin{proof}
$(i)$ Using \eqref{comple}
we have \eq{
    \sum_{j=\max\{k_Q+1,1\}}^\infty\| T_j f\|_{L^p(\Rn \setminus 2Q)}^p 
    &\nonumber\qquad \lesssim \sum_{j=\max\{k_Q+1,1\}}^\infty 2^{j(np-n)}l_Q^{np-n}\\ 
    &\nonumber\qquad \lesssim  2^{k_Q(np-n)}l_Q^{np-n}\\ 
    &\nonumber\qquad \lesssim  l_Q^{-(np-n)}l_Q^{np-n}\\ 
    &\nonumber\qquad\sim 1.
}

$(ii)$ Observe that for $l_Q> 1$, $k_Q \leq 0$ and in this case the statement is trivially true. Hence we assume that $l_Q\leq 1$.

 Now, using \eqref{complementnorm of taylor expansion1} we have 
\begin{align*}
    &\sum_{j=1}^{\max\{k_Q,0\}}\Vert T_j f\Vert^p_{L^p(\Rn \setminus 2Q)} \\ 
    &\nonumber\qquad \lesssim\sum_{j=1}^{\max\{k_Q,0\}} 2^{jm_{\varkappa}(p)p+np/2+jnp\varkappa(1/p-1/2)}2^{jp(N+1)}l_Q^{np-n+p(N+1)}\\ 
    &\nonumber\qquad \lesssim \sum_{j=1}^{\max\{k_Q,0\}} 2^{j(np-n+p(N+1))}l_Q^{np-n+p(N+1)}\\ 
    &\nonumber\qquad \lesssim  2^{k_Q(np-n+p(N+1))}l_Q^{np-n+p(N+1)}\\ 
    &\nonumber\qquad \lesssim  l_Q^{-(np-n+p(N+1))}l_Q^{np-n+p(N+1)}\\ 
    &\nonumber\qquad\sim 1.
\end{align*}
We observe that the same estimates are also valid for the adjoint $T_j^\ast$
since \eqref{comple} and \eqref{complementnorm of taylor expansion1} are valid for the adjoint.
\end{proof}

\begin{Prop}\label{prop:TLp<1}
Suppose that $0<p<1$ and $q\geq p$. Let $a (x,\xi)\in S^{m(p)}_{\rho, \delta}(\Rn)$ be supported outside a neighborhood of the origin and assume that $\varphi\in \textart{F}^k$ {is \emph{SND}}, and satisfies the $L^2$-condition \eqref{eq:L2 condition_old}. Then the \emph{OIO} $T_a^\varphi$ is bounded from $F_{p,q}^{s}(\Rn)$ to $F_{p,q}^s(\Rn)$.
\end{Prop}

\begin{proof}
Observe that it is enough to show the result for $s=0$ and by the inclusion $F^0_{p,p}\xhookrightarrow{}F^0_{p,q}$ it is enough to show that $T_a^\varphi$ is bounded from $F_{p,q}^{0}(\Rn)$ to $F_{p,p}^0(\Rn).$ \\

Compose a Littlewood-Paley piece $\psi_j(D)$ with $T_a^\varphi$ and apply Theorem \ref{thm:left composition with pseudo}. This yields 
\eq{
    \psi_j(D) T_a^\varphi f(x) = T_a^\varphi \psi_j(\nabla_x\varphi(x,D))f(x) +\sum_{\alpha\leq M} 2^{-j\varepsilon}Rf(x),
}
where $\varepsilon>0$ stems from Theorem \ref{thm:left composition with pseudo} and $R$ is an operator of better decay, therefore by Lemma \ref{lem:local and global nonendpoint TL} we have
\eq{
    \norm{\sum_{j=0}^\infty 2^{-j\varepsilon}|Rf(x)|}_{L^p(\Rn)} \sim \|Rf\|_{L^p(\Rn)} \lesssim \| f\|_{F^{-(1-\max(\delta,1/2)-\varepsilon)}_{p,2}(\Rn)}\lesssim \| f\|_{F^0_{p,q}(\Rn)}.
}
So from now on, we will only consider the $F^0_{p,q}\to L^p(\ell^q)$-boundedness of the first term.\\

Let $\mathbf t(\xi)=\nabla_x\varphi(x,\xi)$ and $\eta:\Rn\to\Rn $ be diffeomorphisms such that $\eta(y)\cdot \mathbf t^{-1}(\xi) = y\cdot \xi$. Then we have 
\eq{
    &T_a^\varphi \psi_j(\nabla_x\varphi(x,D))f(x) = \iint_{\Rn\times\Rn} e^{i\varphi(x,\xi)-iy\cdot\xi} \psi_j(\mathbf t(\xi)) a(x,\xi) f(y)\dd y\ddd \xi  
\\&=\frac1{|\det(\nabla t)|}
    \iint_{\Rn\times\Rn} e^{i\varphi(x,\mathbf t^{-1}(\xi))-iy\cdot t^{-1}(\xi)} \psi_j(\xi) a(x,\mathbf t^{-1}(\xi)) f(y)\dd y\ddd \xi  
\\&= 
    \frac{|\det(\nabla \eta)|}{|\det(\nabla t)|}\iint_{\Rn\times\Rn} e^{i\varphi(x,\mathbf t^{-1}(\xi))-iy\cdot \xi} \psi_j(\xi) a(x,\mathbf t^{-1}(\xi)) f(\eta(y))\dd y \ddd \xi 
\\&=:
    T_j (f\circ\eta)(x)
}
Observe that by Theorem \ref{thm:invariance thm} it is enough to consider $T_j f$ from now on.\\

By Lemma \ref{lem:TLatomseq},  $f\in F^0_{p,q}(\Rn)$ has an atomic decomposition 
\eq{
    f(x)=\sum_{\iota=1}^\infty \lambda_\iota  \sum_{ \tilde Q\in \mathcal D_+}b_{\iota,\tilde  Q}\check\Psi^{\tilde Q}(x), \qquad b_{\iota,\tilde Q}\in f^0_{p,q}\cap f^0_{\infty,q},
}
such that
\eq{
    \|f\|_{F^0_{p,q}(\Rn)}\approx\inf_{\set{\lambda_\iota}}\brkt{\sum_{\iota=1}^\infty |\lambda_\iota|^p}^{1/p}.
}

Since $\supp \psi_j\subset \supp \Psi^{\tilde Q}$ it is enough to consider $\tilde Q\in \mathcal D_j$ and hence
\eq{
\norm{\brkt{\sum_{j=1}^\infty|T_jf|^p}^{1/p}}_{L^p(\Rn)} &= \norm{\brkt{\sum_{j=1}^\infty \abs{\sum_{\iota=1}^\infty\lambda_\iota T_j \sum_{\tilde Q\in\mathcal D_j}b_{\iota,\tilde Q}\check\Psi^{\tilde Q}(x)}^p}^{1/p}}_{L^p(\Rn)} \\&\lesssim
\brkt{\sum_{\iota=1}^\infty|\lambda_\iota|^p\int_{\Rn}  \sum_{j=1}^\infty\abs{ T_j \sum_{Q\in\mathcal D_j}b_{\iota,\tilde Q}\check\Psi^{\tilde Q}(x)}^p\dd x}^{1/p} \\&\lesssim
\brkt{\sum_{\iota=1}^\infty|\lambda_\iota|^p}^{1/p}\sup_{\iota\in \mathbb Z_{>0}}\brkt{ \sum_{j=1}^\infty\norm{ T_j \sum_{\tilde Q\in\mathcal D_j}b_{\iota,\tilde Q}\check\Psi^{\tilde Q}(x)}_{L^p(\Rn)}^p}^{1/p}
}

Therefore it is enough to show that one has an expression of the form
\nm{eq:TLgoal}{
    \sum_{j=1}^\infty \| T_j R_{Q,j}\|_{L^p(\Rn)}^p\lesssim 1
}
uniformly in $ Q\in\mathcal D$, where \eq{
    R_{Q,j}(x):=\sum_{\tilde Q\in\mathcal D_j( Q)}b_{\tilde Q}\check\Psi^{\tilde Q}(x).
}
(Recall from Lemma \ref{lem:TLatomseq4} that $R_{Q,j}(x)=0$ for $j< k_{\tilde Q}$.)\\

To show \eqref{eq:TLgoal} we need to show the following three estimates for $0<k<\infty$:
\begin{align}
& \sum_{j=\max\{1,k_Q\}}^\infty\| T_j R_{Q,j}\|_{L^p(2Q)}^p\lesssim 1,\label{eq:TLnewgoal1}\\
& \sum_{j=\max\{1,k_Q\}}^\infty   \|T_j \chi_{2\sqrt n Q} R_{Q,j}\|_{L^p(\Rn \setminus 2Q)}^p\lesssim 1,\label{eq:TLnewgoal2}\\
& \sum_{j=\max\{1,k_Q\}}^\infty\| T_j  \chi_{\Rn \setminus2\sqrt n Q}R_{Q,j}\|_{L^p(\Rn \setminus 2Q)}^p\lesssim 1,\label{eq:TLnewgoal3}
\end{align}

\makeatletter 
\renewcommand{\eqref}[1]{\tagform@{\ref{#1}}}
\def\maketag@@@#1{\hbox{#1}}
\textbf{Step 1 -- Proof of \eqref{eq:TLnewgoal1}}
\makeatother

To see \eqref{eq:TLnewgoal1} consider
\eq{
\| T_j R_{Q,j}\|^p_{L^p(2Q)} &\leq l_Q^{n(2-p)/2p} 2^{2 j m_{\varkappa}(p)} \|2^{-j m_{\varkappa}(p)} T_j R_{Q,j}\|^2_{L^{2}} \\
&\lesssim l_Q^{n(2-p)/2p} 2^{2 j m_{\varkappa}(p)}\| R_{Q,j}\|^2_{L^{2}} \\
&\lesssim l_Q^{n(2-p)/2p} 2^{2k_{Q}  n(1/p-1/2)}2^{2 j m_{\varkappa}(p)}\\
&\lesssim 2^{-k_{Q}n(2-p)/2p} 2^{2k_{Q}  n(1/p-1/2)} 2^{2 j m_{\varkappa}(p)}\\
& = 2^{2 j m_{\varkappa}(p)}
}
for all $0<p<1$, thus we obtain \eqref{eq:TLnewgoal1} by summing in $j$, since $m_{\varkappa}(p)<0$.\\

\textbf{Step 3 -- Proof of \eqref{eq:TLnewgoal2}, \eqref{eq:TLnewgoal3}}\\
\eqref{eq:TLnewgoal2} follows directly from using Lemma \ref{lem:bstarcestimate} and Lemma \ref{lem:TLatomseq4}. \eqref{eq:TLnewgoal3} follows immediately from Lemma \ref{lem:TLatomseq5}. This concludes the proof.
\end{proof}

\section{Transference to Triebel-Lizorkin spaces for \texorpdfstring{$p>2$}{}}\label{subsec:OIO_PRS_transference}

\begin{Def}\label{Influence_set}
In accordance to \emph{Theorem \ref{thm:PRS_proto}}, let $Q$ be a cube and set
    $$\mathcal{E}:= 2Q.$$
\end{Def}
\begin{Lem}\label{OIO_Lemma_PRS}
    Let $n\geq 1$, $0<p<1$. Let $0\leq \rho\leq 1$ and $0\leq \delta< 1$. Assume that $\varphi\in \textart F^k$ is \emph{SND} and satisfies the $L^2$-condition \eqref{eq:L2 condition_old} for $0<k<\infty$. 
Let $a(x,\xi)\in S^{m_{\varkappa}(p)}_{\rho,\delta}(\Rl^n)$, and 
\begin{equation*}
    K_{j}(x,y) = \int_{\Rl^n} e^{i(\varphi(x,\xi)-y\cdot \xi)}\,\sigma_{j}(x,\xi)  \ddd \xi.
\end{equation*} 
be the kernel of $T_j$. Then for all $\eps>0$ we have for $\at$ a $h^p$-atom supported in $Q$ that
    \begin{enumerate}
        \item[$(i)$] If $2^{-j}\leq l_Q\leq 1$, then
        \begin{equation}\label{small_k_lift_equation}
            \int_{\Rn} |T_j^*\Psi_j(D) \at| \dd x
    \lesssim 2^{j m_{\varkappa}(p) \frac{1}{\frac{2}{p}-1}}l_Q^{-\eps}+l_Q^{\frac{1}{\frac{2}{p}-1}(n-n/p)}2^{\frac{1}{\frac{2}{p}-1}j(n-n/p)},
        \end{equation}
        \item[$(ii)$] If $l_Q\geq 1$, then
         \begin{equation}\label{small_k_lift_equation2}
            \int_{\Rn} |T_j^*\Psi_j(D) \at| \dd x
            \lesssim (2^{j m_{\varkappa}(p)}+l_Q^{n-n/p}2^{j(n-n/p)})^{\frac{1}{\frac{2}{p}-1}}.
        \end{equation}
    \end{enumerate}
\end{Lem}

\begin{proof}
We begin by proving \eqref{small_k_lift_equation}, to this end let $w$ be a real number such that $1/w=1/p-1/2$, and split $\Rn = \mathcal E \cup \Rn\setminus \mathcal E$. Then
\begin{align*}
    \Big(\int_{\mathcal E}  |T_j^*\Psi_j(D) \at|^p \dd x\Big)^{1/p} 
    &\lesssim \| 1 \|_{L^w(\mathcal E)} \Big(\int_{\Rn}  |T_j^*\Psi_j(D) \at|^2 \dd x\Big)^{1/2} \\
    &\lesssim |Q|^{1/p-1/2}2^{jm_\varkappa(p)}\Big(\int_{\Rn}  | \at(x)|^2 \dd x\Big)^{1/2}\lesssim 2^{jm_{\varkappa}(p)}.
\end{align*}
Now, recall that  
$$g(x)=\jap{2^{j\varkappa}|x-\bar y|}^{-M}$$
and that we have for $w\geq 1$,
\begin{equation}\label{metric norm}
    \|g\|_{L^w(\Rl^n)}\lesssim 2^{nj\varkappa/w}.
\end{equation}

Next H\"older's and Minkowski's inequalities and Lemma \ref{L2 of kernel lemma} yield that
\begin{align*}
    &\Big(\int_{\Rn\setminus\mathcal E}  |T_j^*\Psi_j(D) \at|^p \dd y\Big)^{1/p} \lesssim \Big\| \frac{1}{g(x)} \int_B  \overline{ K_j(y,x)}\Psi_j(D)\at(y) \dd y\Big\|_{L^2(\Rl^n)}\|g\|_{L^w(\Rl^n)}\\
    &\nonumber \lesssim  \int_B\Big\| \frac{1}{g(x)} \overline{ K_j(y,x)}\Big\|_{L^2(\Rl^n)} |\Psi_j(D)\at(y)| \dd y\,\|g\|_{L^w(\Rl^n)}\\
    &\nonumber \lesssim  2^{-nj\varkappa/w} \int_B\Big\| \frac{1}{g(x)}\overline{ K_j(y,x)}\Big\|_{L^2(\Rl^n)} |\Psi_j(D)\at(y)| \dd y\\
    &\nonumber \lesssim  l_Q^{n-n/p}2^{j(m_\varkappa(p)+n/2-n\varkappa(1/p-1/2))}\\
    &\nonumber =  l_Q^{n-n/p}2^{j(n-n/p)}
\end{align*}
where we have also use that $\at$ is an $h^p$-atom and that $\Psi_j(D)$ is $L^2$-bounded uniformly in $j$.\\

\textbf{Interpolation step:}\\

Using the $L^2$ result (Theorem \ref{calderon-vaillancourt lemma}) we obtain by Riesz-Thorin interpolation and $\frac{1}{p_t}=\frac{t}{2}+ \frac{1-t}{p}$ with $0<p<1$ and $0<k< 1$ that
\begin{align*}
    \Big(\int_{\Rn} |T_j^*\Psi_j(D) \at|^{
    p_t} \dd x\Big)^{1/p_t}
    &\lesssim (2^{j m}+l_Q^{n-n/p}2^{j(n-n/p)})^{(1-t)}\\
    &\lesssim 2^{j m (1-t)}+l_Q^{(1-t)(n-n/p)}2^{(1-t)j(n-n/p)}.
\end{align*}

Taking $p_t=1$ and $t=\frac{1-\frac{1}{p}}{\frac{1}{2}-\frac{1}{p}}$ we have $1-t=\frac{1}{\frac{2}{p}-1}$ and therefore
\begin{align*}
    \int_{\Rn} |T_j^*\Psi_j(D) \at| \dd x
    &\lesssim (2^{j m}+l_Q^{n-n/p}2^{j(n-n/p)})^{\frac{1}{\frac{2}{p}-1}}.
\end{align*}

Thus for $l_Q>1$ we have
\begin{align*}
    \int_{\Rn} |T_j^*\Psi_j(D) \at| \dd x
    &\lesssim 2^{j m \frac{1}{\frac{2}{p}-1}}+2^{\frac{1}{\frac{2}{p}-1}j(n-n/p)}
\end{align*}

and for $l_Q<1$ and $\eps>0$ we have
\begin{align*}
    \int_{\Rn} |T_j^*\Psi_j(D) \at| \dd x
    &\lesssim 2^{j m \frac{1}{\frac{2}{p}-1}}l_Q^{-\eps}+l_Q^{\frac{1}{\frac{2}{p}-1}(n-n/p)}2^{\frac{1}{\frac{2}{p}-1}j(n-n/p)}.
\end{align*}
\end{proof}
 \begin{Th}\label{thm:PRS}
 Let $n\geq 1$. Assume that $\varphi\in \textart F^k$ is \emph{SND} and satisfies the $L^2$-condition \eqref{eq:L2 condition_old} for $0<k<\infty$. 
Then if $a(x,\xi)\in S^{m_{\varkappa}(p)}_{\rho,\delta}(\Rl^n)$, and $T_j$ are as in \emph{Lemma \ref{OIO_Lemma_PRS}}. Finally, we let $2 < p$, $0<q\leq \infty$, and $b>0$.  Assume that the operators $T_j$ satisfy
\begin{align}
&\sup_{
j>0}
2^{jb/p}\| T_j\|_{L^p\to L^p} \lesssim 1,\label{eq:PRS1}
\\&
\sup_{
j>0}
2^{jb/2}\|T_j\|_{L^2\to L^2} \lesssim 1.\label{eq:PRS2}
\end{align}

Then, the following inequality holds:
\eq{
\norm{ \brkt{\sum_{j=0}^\infty 
2^{jbq/p}|\Psi_j(D) T_jf_j|^q}^{1/q}}_{L^p(\Rn)} \lesssim \brkt{ \sum_{j=0}^\infty \|f_j\|_{L^p(\Rn)}^p}^{1/p}.
}
where $\Psi\in \mathscr{S}(\Rl^n)$, $\Psi_j(D):=\Psi(2^{-j}D)$ and $f_j$ is a sequence of functions.

\end{Th}

\begin{proof}
We consider a measurable set $\mathcal E$ (as defined in Definition \ref{Influence_set}). We have the following inequality:
\begin{equation*}
|\mathcal E| \lesssim l_Q^n,
\end{equation*}
which means that all the hypotheses of Theorem \eqref{thm:PRS_proto} are satisfied, with $q=2$. Note that the condition $b<n$ in Theorem \eqref{thm:PRS_proto} is not necessary here, because it stems from Lemma 2.2 in \cite{PRS}, which we substitute with our Lemma \eqref{OIO_Lemma_PRS} (and duality). We also do not need the assumption \eqref{eq:PRS} about the kernel, because we show Lemma \ref{OIO_Lemma_PRS} using a different argument. Therefore, we can obtain the result by applying the same method as in \cite{PRS} to the adjoint $T_j^*$, and then interpolating with $L^2(\Rn)$.
\end{proof}
\section{Triebel-Lizorkin estimates}\label{main TL estim section}

In this section we use the lifting results from section \ref{subsec:OIO_lowP_transference} and \ref{subsec:OIO_PRS_transference} to lift Theorem \ref{thm:hp-Lp_oio} to Triebel-Lizorkin boundedness.

\begin{figure}\label{Firgure triebel lizorkin}
\begin{tikzpicture}

    \draw[thick,->] (0,0) -- (7,0) node[anchor=north west] {$p$};
    \draw[thick,->] (0,0) -- (0,7) node[anchor=south east] {$q$};

    \draw (0,0) -- (7,7);
    \draw[dashed] (2,7) -- (2,0);
    \draw[dashed] (0,2) -- (7,2);
    \draw[dashed] (4,0) -- (4,7);
    \draw[dashed] (0,4) -- (7,4);
    
    \draw[thick] (2 cm,2pt) -- (2 cm,-2pt) node[anchor=north] {$2$};
    \draw[thick] (4 cm,2pt) -- (4 cm,-2pt) node[anchor=north] {$1$};
    \draw[thick] (2pt,2 cm) -- (-2pt,2 cm) node[anchor=east] {$2$};
    \draw[thick] (2pt,4 cm) -- (-2pt,4 cm) node[anchor=east] {$1$};
    \draw (0,0) node[anchor=north east] {$\infty$};

    \foreach \x in {2.1,2.25,...,4}{
        \draw[blue] (\x cm,0) -- (\x cm,\x cm);
    }
    \foreach \x in {4.1,4.25,...,6.90}{
        \draw[red] (\x cm,0) -- (\x cm,\x cm);
    }
    \foreach \y in {0,0.15,...,2}{
        \draw[blue] (0,\y cm) -- (\y cm,\y cm);
    }
    \foreach \y in {2,2.15,...,6.8}{
        \draw[blue] (0,\y cm) -- (2 cm,\y cm);
    }

\end{tikzpicture}
\caption{Triebel-Lizorkin boundedness for OIO:s with their respective critical decay. The blue horizontal lines illustrate Theorem \ref{thm:PRS}. The blue vertical lines illustrate the boundedness results obtained by applying a duality argument for the case $p>2$. The red vertical lines illustrate Proposition \ref{prop:TLp<1} together with interpolation with the vertical blue area.}\label{pic:TLendpointresults}
\end{figure}
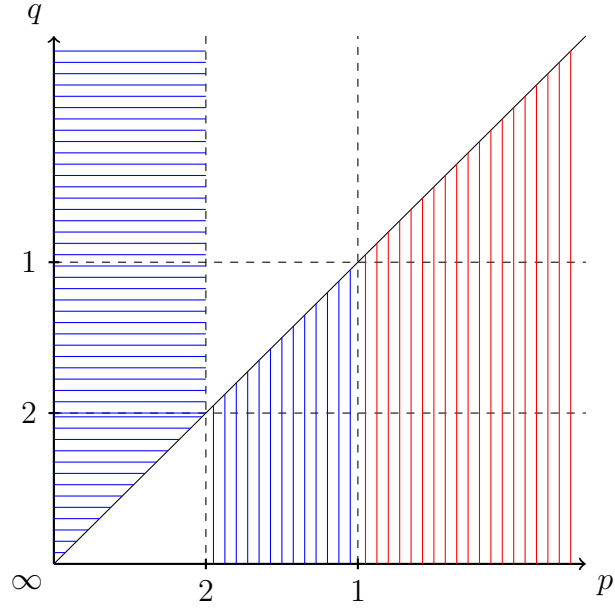

\begin{Th}\label{thm:TL_tdphase_k}
    Let $0\leq \rho\leq 1$, $0\leq\delta<1$, $k>0$, $0<\mu\leq 1$, and $0<q\leq \infty$.
     Assume furthermore that $\varphi\in \textart F^k$ {is \emph{SND}}, and satisfies the \emph{LF}$(\mu)$ condition \eqref{eq:LFmu} and the $L^2$-condition \eqref{eq:L2 condition_old}. For $a\in S^{m(p)}_{\rho, \delta}(\Rl^n)$, and let $T_a^\varphi$  be the associated \emph{OIO}.
    If $s\in \Rl$ and either one of the following cases holds
    \begin{align*}
    (i)\quad &2<p<\infty \text{  when }0<q\leq p,\\
    (ii) \quad &\frac{n}{n+\mu}<p<2 \text{ when }p\leq q,\\
    (iii)\quad &p=q=2,
    \end{align*}
    then the \emph{OIO} $T_a^\varphi$ is bounded from $F_{p,q}^{s}(\Rn)$ to $F_{p,q}^s(\Rn)$
\end{Th}

\begin{proof}
We separate the operator into a low and a high-frequency part. The result for the low-frequency part follows from Theorem \ref{thm:low_freq_TL_BL_OIO}, so we only consider the high-frequency part from now on. \\

Observe also that the contents of $(iii)$ is contained in Theorem \ref{thm:Calderon-Vaillancourt for OIOs}. So from now on we only need to show $(i)$ and $(ii)$.\\

We split the proof into different ranges of $p$ and $q$, the two parts of the proof correspond to the blue and the red regions in {\bf Figure \ref{pic:TLendpointresults}}, respectively.\\

\textbf{Part 1 -- Proof when $\boldsymbol{p > 2}$ and $\boldsymbol{p\geq q>0}$}\\
We use Theorem \ref{thm:left composition with pseudo} to write 
\eq{
    \psi_j(D) T_a^\varphi f(x) = T_a^\varphi \psi_j(\nabla_x\varphi(x,D))f(x) + Rf(x),
}
The operator $R$ is bounded by Lemma \ref{lem:local and global nonendpoint TL}. So from now on we will only consider the first term. Denote $T_j:= T_a^\varphi \psi_j(\nabla_x\varphi(x,D)).$\\

Hence we can use Theorem \ref{thm:PRS} with $b:=n(1-\varkappa)$ and $S_j:= T_j(1-\Delta)^{-b/2p}$ to prove the desired result.
Observe that the $h^p\to L^p$ boundedness (Proposition \ref{thm:hp-Lp_oio}) and the $L^2$ boundedness (Theorem \ref{thm:Calderon-Vaillancourt for OIOs})of $T_a^\varphi$ yield \eqref{eq:PRS1} and \eqref{eq:PRS2} respectively.\\

Theorem \ref{thm:PRS} immediately yields that
\eq{
    &\norm{ \brkt{\sum_{j=0}^\infty 
    2^{jbq/p}|\Psi_j(D) T_jf_j|^q}^{1/q}}_{L^p(\Rn)}=\norm{ \brkt{\sum_{j=0}^\infty 
    2^{jbq/p}|\Psi_j(D) S_jf_j|^q}^{1/q}}_{L^p(\Rn)} 
    \\&
    \lesssim \brkt{ \sum_{j=0}^\infty \|f_j\|_{L^p(\Rn)}^p}^{1/p}.
}

Thus $S_j:B^0_{p,p}\to F^{b/p}_{p,q}$ which immediately implies that $T_j:B^{b/p}_{p,p}\to F^{b/p}_{p,q}$. Now the assertion follows from the facts that $F^0_{p,q}\xhookrightarrow{} F^0_{p,p}=B^0_{p,p}$ and the calculus using Bessel potentials and Theorem \ref{thm:left composition with pseudo}.\\

\textbf{Part 2 -- Proof when $\boldsymbol{1<p<2}$ and $\boldsymbol{p\leq q}$}\\
Using that the operator is self-adjoint we can also obtain Theorem \eqref{thm:PRS} for the adjoint, then apply the arguments from part 1.\\
 
\textbf{Part 3 -- Proof when $\boldsymbol{0<p<1, p\leq q\leq\infty}$}\\
By Proposition \ref{prop:TLp<1} one obtains the result for $0<p<1$ and $q\geq p$.\\

\textbf{Part 3.1 -- Proof when $\boldsymbol{0<p\leq q<1}$}\\
In this case using Riesz-Thorin interpolation with $F^s_{p,p}=B^s_{p,p},$ yields the result. \\

\textbf{Part 3.2 -- Proof when $\boldsymbol{0<p<1\leq q\leq \infty}$}\\
Here Riesz-Thorin interpolation with {\bf{Part 2}} above, yields the result. \\

Now notice that $(1-\Delta)^{\frac s2} T^{\varphi}_a(1-\Delta)^{-\frac s2}$ is a similar operator associated to an amplitude in $S^{m_c(p)}(\Rn)$ and phase $\varphi$, and hence bounded from $F^0_{p,q}(\Rn)$ to itself. Therefore using the fact that the operator $(1-\Delta)^{\frac{s}{2}}$ is an isomorphism from $F^s_{p,q}(\Rn)$ to $F^0_{p,q}(\Rn)$ for $0<p\leq \infty$, we obtain the desired result. 
\end{proof}

\subsection{Triebel-Lizorkin estimates related to forbidden amplitudes}\label{subsec:OIO_Forbidden}
\quad\\

\noindent It is well known that the oscillatory integral operators with amplitudes in $S^{m}_{1,1}(\Rl^n)$ fail to be $L^2$-bounded. However, one may show that these operators are Sobolev-bounded for $H^s(\Rn)$ for $s>0$ (see \cite{CIS}), and the pseudodifferential case goes back to E. Stein and independently by Y. Meyer. In this section, we state and prove two results about the boundedness of oscillatory integral operators in Triebel-Lizorkin spaces with amplitudes in $S^{m}_{1,1}(\Rl^n)$.\\

Now we turn to the boundedness of OIOs with forbidden amplitudes in the class $S^m_{1,1}(\mathbb{R}^n)$ and ask whether they are bounded on Triebel-Lizorkin spaces. 
 
\begin{Th}\label{thm:Sobolev_oio}
Let $n\geq 1$, $k>0$, $s>0$, $0<\mu\leq 1$. Assume that \linebreak$a(x,\xi)\in S^{-kn|1/p-1/2|}_{1, 1}(\Rl^n)$, $\varphi\in \textart F^k$ {is \emph{SND}}, satisfies the \emph{LF}$(\mu)$-condition and the conditions in \eqref{eq:composition_conditions}, 
 and the $L^2$-condition \eqref{eq:L2 condition_old}.  
If $s>n\big(\frac{1}{\min\{1, p, q\}}-1\big)$,
    and either one of the following cases hold 
    \begin{enumerate}
        \item[$(i)$] $2<p<\infty$ when $0<q\leq p,$\\
        \item[$(ii)$] $\frac{n}{n+\mu}<p<2$ when $p\leq q,$\\
        \item[$(iii)$] $p=q=2,$
    \end{enumerate}
    or if $s>\frac{n}{p}$ with $q=\infty$, then the \emph{OIO} $T_a^\varphi$ is bounded from $F^s_{p,q}(\Rl^n)\to F^s_{p,q}(\Rl^n).$ \\
\end{Th}

\begin{proof}
The proof is similar to the proof of Theorem 5.13 in \cite{ATW}, and therefore we shall only highlight the differences here. The only differences appear in connection to (100). At this point in the proof one instead uses Theorem \ref{thm:TL_tdphase_k} instead of Theorem 5.15. The rest of the proof remains the same.
\end{proof}

\section{Besov-Lipschitz estimates}\label{main results in BL}

In this section we include both the Besov-Lipschitz boundedness results of OIO's with amplitudes in $S^m_{\rho,\delta}$ for all $0\leq \rho\leq 1$ and $0\leq \delta\leq 1$. We begin with the classical amplitudes were $\delta<1$.

\begin{Th}\label{thm:BL_tdphase_k_less_than_one}
    Let $0\leq \rho\leq 1$ and $0\leq\delta<1$, $k>0$, $0<\mu\leq 1,$ $0<q\leq\infty$.
    Assume that $\varphi\in \textart F^k$ {is \emph{SND}}, satisfies the \emph{LF}$(\mu)$-condition and the conditions in \eqref{eq:composition_conditions}, 
    and the $L^2$-condition \eqref{eq:L2 condition_old}, and $a\in S^{m(p)}_{\rho, \delta}(\Rl^n)$, and let $T_a^\varphi$  be the associated \emph{OIO}.
    Then $T_a^\varphi$ is a bounded operator from $B^s_{p,q}(\Rl^n)$ to $B^s_{p,q}(\Rl^n)$ for all $s\in\Rl$ for $\frac{n}{n+\mu}<p<\infty$.
\end{Th}

\begin{proof}
    Let $\chi\in C_c^\infty(\Rn)$ be supported in $\{\xi:|\xi|\lesssim 1\}$, and write
    \begin{align*}
        T_a^\varphi f(x) &: =
        \int_{\Rn} e^{i\varphi(x,\xi)} a(x,\xi) \widehat{f}(\xi) (1-\chi(\xi)) \ddd\xi 
        +
        \int_{\Rn} e^{i\varphi(x,\xi)} a(x,\xi) \widehat{f}(\xi) \chi(\xi) \ddd\xi\\
        &=T_{\text{high}} f(x)+T_{\text{low}} f(x).
    \end{align*}
    The boundedness of $T_{\text{low}}$ follows from the low frequency result Theorem \ref{thm:low_freq_TL_BL_OIO}. The boundedness of $T_{\text{high}}$ follows immediately by the Besov-Lipschitz lift Theorem \ref{theorem:BL-booster theorem} and Proposition \ref{thm:hp-Lp_oio}.
\end{proof}

In this next result we obtain an estimate relating to OIO's with amplitudes in $L^{\infty}S^{m}_{0}(\Rl^n)$ (see Definition \ref{symbol class Sm no regularity}), observe that this class of amplitudes contains the forbidden H\"ormander class amplitudes $S^{m}_{0,1}(\Rl^n)$ and $S^{m}_{1,1}(\Rl^n)$.

\begin{Th}\label{thm:Sobolev_oio0}
Let $n\geq 1$, $s\in \Rl$, $1<p< \infty$, $ 0< q\leq \infty$. Assume furthermore that $\varphi(x, \xi)\in\textart F^k$ is an \emph{SND} phase of order $k>0$, $a\in L^{\infty}S^{m}_{0}(\Rl^n)$ with $$m< n\Big(kp-1+\Big(1-\frac{1}{p}\Big)w(k,\rho)\Big),$$ and
the \emph{OIO} $T_a^\varphi$ is bounded from $B^s_{p,q}(\Rl^n)\to B^s_{p,q}(\Rl^n).$ $($Where $w(k,\rho)$ is as defined in \emph{Lemma \ref{kernel lemma})}
\end{Th}

\begin{proof}

We separate the operator into a low and a high-frequency part. The result for the low-frequency part follows from Theorem \ref{thm:low_freq_TL_BL_OIO}, so we only consider the high-frequency part from now on. \\

Following the proof of Theorem \ref{thm:hp-Lp_oio} we obtain the kernels $K_j^\nu$ with the associated kernel estimate \eqref{main kernel estimate for AAW 1}.\\

Let $T_j^\nu$ be the operators corresponding to the kernels $K_j^\nu$. One observes that
\begin{align*}
 |T_j^\nu f_j(x)|^{r}&=\Big|\int_{\Rn} K_j^\nu(x,y)f_j(y)\dd y\Big|^{r}\\
&=\Big|\int_{\Rn} K_j^\nu(x,y)\sigma(\nabla_\xi \varphi(x,\xi_j^\nu)-y)\frac{1}{\sigma(\nabla_\xi \varphi(x,\xi_j^\nu)-y)}f_j(y)\dd y\Big|^{r},
\end{align*}
with weight functions $\sigma(y)$ which will be chosen
momentarily. Therefore, H\"older's inequality with $\frac{1}{r}+ \frac{1}{r'}=1$, $r,r'>1$ yields
\begin{align}
    &|T_j^\nu f_j(x)|^{r}\leq
\Big( \int_{\Rn} |K_j^\nu(x,y)|^{r'} | \sigma(\nabla_\xi \varphi(x,\xi_j^\nu)-y)|^{r'} \dd y\Big)^{\frac{r}{r'}}\label{eq2.14}\\
  &\qquad\times\Big( \int_{\Rn} \frac{| f_j(y)|^{r}}{|\sigma(\nabla_\xi \varphi(x,\xi_j^\nu)-y)|^{r}}\dd y\Big).\nonumber
\end{align}
where $\sigma$ is defined by
\begin{equation*}
\sigma(y)=\begin{cases}
1, &| y| \leq 1; \\
| y|^{\lambda}, & | y| >1.
\end{cases}
\end{equation*}
Observe further that
\begin{align*}
    &\int_{\Rn} |K_j^\nu(x,y)|^{r'} | \sigma(\nabla_\xi \varphi(x,\xi_j^\nu)-y)|^{r'} \dd y\\
    &\lesssim \int_{\Rn} \frac {2^{r'j(m )} 2^{r'jn(1-k)}}{\jap{2^{jw(k,\rho)}(\nabla_\xi \varphi(x,\xi_j^\nu)-y)}^{Mr'}} | \sigma(\nabla_\xi \varphi(x,\xi_j^\nu)-y)|^{r'} \dd y\\
    &= \int_{\Rn} \frac {2^{r'j(m )} 2^{r'jn(1-k)}}{\jap{2^{jw(k,\rho)}y}^{Mr'}} | \sigma(y)|^{r'} \dd y\\
    &\lesssim  \int_{|y|> 1} \frac{2^{jr'm} 2^{-(M-\lambda) r'w(k,\rho)}2^{jnr'(1-k)}}{|y|^{(M-\lambda) r'}} \dd y +\int_{|y|\leq 1} \frac{2^{jr'm} 2^{jnr'(1-k)}}{\jap{2^{jw(k,\rho)}y}^{Mr}} \dd y\\
    &\lesssim  2^{jr'm} 2^{-(M-\lambda) r'w(k,\rho)}2^{jnr'(1-k)} +2^{jr'm} 2^{jn(r'(1-k)-w(k,\rho))}.
\end{align*}

Furthermore, using (17) in \cite[p. 57]{Stein}, we have
\begin{equation*}
\int_{\Rn} \frac{| f_j(y)|^{r}\dd y}{|\sigma(\nabla_\xi \varphi(x,\xi_j^\nu)-y)|^{r}}\lesssim \big(\mathcal{M}_{r} f_j(\nabla_\xi \varphi(x,\xi_j^\nu))\big)^{r}
\end{equation*}

with a constant that only depends on the dimension $n$. Thus
\eqref{eq2.14} yields
\begin{equation*}
    |T_j f_j(x)| \lesssim 2^{jr'm/r} 2^{jn(r'(1-k)-w(k,\rho))/r} \sum_{\nu} \mathcal{M}_{r} f_j(\nabla_\xi \varphi(x,\xi_j^\nu)).
\end{equation*}
From which it now follows for $r<p$, the SND condition and the Hardy-Littlewood maximal theorem that
\begin{align*}
    \|T_j f_j\|_{L^p(\Rn)} &\lesssim 2^{jr'm/r} 2^{jn(r'(1-k)-w(k,\rho))/r} \sum_{\nu} \big\|\mathcal{M}_{r} f_j(\nabla_\xi \varphi(\cdot,\xi_j^\nu))\big\|_{L^p(\Rn)}\\
    &\lesssim 2^{jr'm/r} 2^{jn(r'(1-k)-w(k,\rho))/r} 2^{nj k}\|f_j\|_{L^p(\Rn)},
\end{align*}
and therefore, by Minkowski's integral inequality,
\begin{align*}
    \|T_a^\varphi f\|_{B^s_{p,q}}&\lesssim \Big(\sum_{j=1}^\infty 2^{jqs}\|T_j f_j\|^q_{L^p(\Rn)}\Big)^{1/q}\\
    &\lesssim \Big(\sum_{j=1}^\infty 2^{jqs}2^{qjr'm/r} 2^{qjn(r'(1-k)-w(k,\rho))/r}2^{qnjk}\|f_j\|^q_{L^p(\Rn)}\Big)^{1/q}\\
    &\lesssim \|f\|_{B^s_{p,q}(\Rn)},
\end{align*}
whenever
$$m< n\Big(kp-1+\Big(1-\frac{1}{p}\Big)w(k,\rho)\Big).$$
\end{proof}

 \bibliographystyle{siam}
 \bibliography{references}

\end{document}